\documentclass[final,onefignum,onetabnum]{siamart190516}

\usepackage{verbatim}
\usepackage{lipsum}
\usepackage{amsfonts}
\usepackage{graphicx}
\usepackage{subfigure}
\usepackage{epstopdf}
\usepackage{algorithmic}
\usepackage{amsmath,amssymb}
\ifpdf
  \DeclareGraphicsExtensions{.eps,.pdf,.png,.jpg}
\else
  \DeclareGraphicsExtensions{.eps}
\fi
\renewcommand\theequation{\thesection.\arabic{equation}}
\theoremstyle{plain}

\newtheorem{remark}{\textbf{Remark}}[section]

\newcommand{\eps}{\epsilon}

\newcommand{\bm}{\boldsymbol}

\newcommand{\Grad}[1]{\nabla #1}

\newcommand{\be}{\begin{equation}}
\newcommand{\ee}{\end{equation}}

\newcommand{\bse}{\begin{subequations}}
\newcommand{\ese}{\end{subequations}}
\def\benl{\begin{eqnarray*}}
\def\eenl{\end{eqnarray*}}

\def\be{\bm{e}}

\def\bx{\bm{x}}

\def\bz{\bm{z}}

\def\bmu1{\bm{\mu_1}}

\def\ln{{\rm ln}}

\newcommand{\ben}{\begin{eqnarray}}
\newcommand{\een}{\end{eqnarray}}
\newcommand{\beq}{\begin{equation}}
\newcommand{\eeq}{\end{equation}}
\newcommand{\bea}{\begin{array}}
\newcommand{\eea}{\end{array}}
\newcommand{\bef}{\begin{figure}(H)_h}
\newcommand{\eef}{\end{figure}}

\newsiamremark{hypothesis}{Hypothesis}
\crefname{hypothesis}{Hypothesis}{Hypotheses}
\newsiamthm{claim}{Claim}

\title{A new Lagrange multiplier approach for constructing structure-preserving schemes, II. bound preserving}

\author{Qing Cheng \thanks{Department of Mathematics, Purdue University, West Lafayette, IN 47907, USA  (cheng573$@$purdue.edu).}
\and  Jie Shen\thanks{Department of Mathematics, Purdue University, West Lafayette, IN 47907, USA  (shen7$@$purdue.edu)}
}
  
\usepackage{amsopn}

\begin{document}
\bibliographystyle{plain}
\graphicspath{ {Figures_porous/} }
\maketitle

\begin{abstract}
In the second part of this series, we use the  Lagrange multiplier approach proposed in the first part \cite{CheS21} to construct efficient and accurate bound and/or mass preserving schemes for  a class of semi-linear and quasi-linear parabolic equations.   
We establish    stability results under a general setting, and carry out an error analysis for a second-order bound preserving scheme with a hybrid spectral discretization in space. We apply our approach to several typical PDEs which preserve bound and/or mass, also present ample numerical results to validate our approach.

\end{abstract}

\begin{keywords}
bound preserving;  mass conservation; KKT conditions;  Lagrange multiplier; stability; error analysis
\end{keywords}

\begin{AMS}
65M70; 65K15; 65N22
\end{AMS}

\section{Introduction}

Solutions of partial differential equations (PDEs) arising from sciences and engineering applications are often required to be positive or to remain in a bounded interval.   It is beneficial, and often necessary,  that their numerical approximations preserve the positivity or bound at the discrete level.
In recent years, a large effort has been devoted to construct  bound preserving schemes for various problems. 

In the first part of this series \cite{CheS21}, we constructed a class of positivity preserving schemes using a new Lagrange multiplier approach. A main objective of this paper is to extend the approach in  \cite{CheS21} to construct bound preserving schemes for a class of nonlinear PDEs in the following form: 
\begin{equation}\label{strong}
 u_t + \mathcal{L}  u+\mathcal{N}(u)=0,
\end{equation}
with  suitable initial and boundary conditions, where $\mathcal{L}$ is a linear or nonlinear  non-negative operator and $\mathcal{N}(u)$  is a semi-linear or quasi-linear operator. 
We assume that the solution of 
 \eqref{strong} is bound preserving, i.e., $a\le u(\bx,0)\le b$ for all $\bx\in \Omega$, then $a\le u(\bx,t)\le b$ for all $(\bx,t)\in \Omega\times (0,T)$. 
 
 There are a large body of work devoted to construct positivity/bound preserving schemes for \eqref{strong}. We refer to  the first part of this series \cite{CheS21} (and the references therein) for a summary of existing approaches for constructing positivity/bound preserving schemes. In particular, large efforts have been devoted to construct spatial discretization for \eqref{strong} such that  the resulting  numerical scheme satisfies a discrete maximum principle (cf., for instance,   \cite{du2019maximum,ciarlet1970discrete,ciarlet1973maximum,droniou2011construction,liu2014maximum,MR3913649,liao2020second,MR4166061}, and the review paper in \cite{MR4253790} for a up-to-date summary in this regard).

Given a  generic spatial discretization of \eqref{strong}:
\begin{equation}\label{appro_strong}
\begin{split}
\partial_t u_h + \mathcal{L}_h  u_h +\mathcal{N}_h(u_h)=0,
\end{split}
\end{equation}
where $u_h$ is in certain finite dimensional approximation space $X_h$  and $\mathcal{L}_h$ is a certain approximation of $\mathcal{L}$. In general, the solution $u_h$,  if it exists, may not be bound preserving. Oftentimes, \eqref{appro_strong} may not be well posed if the values of $u_h$  go outside of $[a,b]$. For example, a direct finite elements or spectral approximation to the  Allen-Cahn or Cahn-Hilliard equation with logarithmic potential may not be well posed. Instead of using special  spatial discretizations which  satisfy a discrete maximum principle,
 we aim to develop a bound preserving approach which can be used for a large class of spatial discretizations. To preserve positivity, it suffices to introduce a Lagrange multiplier $\lambda_h$. But to preserve bound, we need to introduce an additional quadratic  function $g(u)=(b-u)(u-a)$, and consider the
 following expanded system with a Lagrange multiplier $\lambda_h$:
\begin{equation}\label{kkt}
\begin{split}
&\partial_t u_h +  \mathcal{L}_h  u_h+\mathcal{N}_h(u_h)=\lambda_hg'(u_h),\\
&\lambda_h \ge 0,\;g(u_h)\ge 0,\;\lambda_h g(u_h)=0.
\end{split}
\end{equation}
The second equation in \eqref{kkt} represents the well-known KKT conditions \cite{ito2008lagrange,facchinei2007finite,harker1990finite,bergounioux1999primal} for constrained minimization.  
The problem \eqref{kkt} can be viewed as an approximation to \eqref{strong}, it can also be viewed as a discrete problem without a background PDE, e.g., coming from a discrete constrained minimization problem.

Existing approaches for \eqref{kkt} usually start with  an implicit  time discretization scheme so that the nonlinear system at each time step can still be interpreted as a  constrained minimization, then  apply a suitable iterative procedure (cf. \cite{van2019positivity}). As in \cite{CheS21}, we shall use a different approach which decouples the computation of  Lagrange multiplier $\lambda_h$ from that of $u_h$, leading to a much more efficient algorithm.

We recall that for positivity preserving, we simply use $g(u_h)=u_h$ in the above formulation.  However, for bound preserving, the nonlinear nature of $g(u_h)$ makes it much harder to prove  stability in norms involving derivatives, and mass conservation whenever is necessary. On the other hand, since the numerical solutions remain to be bounded by construction, this allows us to derive more precise stability results, which in turns enable us to obtain optimal error estimates for both semi-linear and quasi-linear PDEs. More precisely, the 
bound preserving schemes that we construct based on the operator splitting approach enjoy all  advantages of the positivity preserving schemes in \cite{CheS21}, and furthermore, thanks to the bound preserving property, it allows us to prove a more precise stability result (see Theorem 3.1) and to establish rigorous error estimates for a class of semi-linear and quasi-linear dissipative equations (see Theorem 4.1).

We would like to point out that the schemes  constructed in this paper include the usual cut-off approach \cite{MR3062022}  as a special case. Therefore, our presentation provides an alternative interpretation of the   cur-off approach, and allows us to construct  new cut-off  implicit-explicit (IMEX) schemes with mass conservation. 

To validate our schemes, we  apply our new schemes to a variety of  problems with bound preserving solutions, including the Allen-Cahn \cite{allen1979microscopic} and Cahn-Hilliard \cite{cahn1958free} equations and a class of Fokker-Planck equations \cite{risken1996fokker}.

The remainder of the paper is organized as follows.  In Section 2, we construct bound  preserving schemes for general nonlinear systems \eqref{strong} using the  Lagrange multiplier approach. For problems which also conserve mass, we modify our bound preserving schemes so that they also conserve mass. In Section 3, we  restrict ourselves to second-order parabolic type equations, and  establish  a  stability result for, as an example,  the second-order    scheme with mass conservation. In Section 4,  we consider a hybrid  spectral method as an example to carry out an error analysis for a fully discretized second-order  scheme. In Section 5, we describe applications of our schemes to several typical PDEs with bound and/or mass preserving properties. In Section 6, we present some numerical simulations  to validate the accuracy and robustness of  our schemes. And we conclude with some  remarks   in the  final section.

\section{Bound-preserving schemes}
We construct in this section efficient bound preserving schemes for solving   \eqref{kkt}. The key is to adopt an operator splitting approach in which a standard scheme, which is not bound preserving,  is used in the first step, while in the second step, the solution is made bound preserving with a simple yet consistent procedure.

We shall first describe a generic spatial discretization with nodal Lagrangian basis functions, followed by time discretization without and with mass conservation.


 Let $\Sigma_h$ be a set of mesh points or collocation points in $\bar\Omega$. Note that $\Sigma_h$ should not include the points at the part of the boundary where a Dirichlet (or essential) boundary condition is prescribed, while it should include the points at the part of the boundary where a Neumann or mixed  (or non-essential) boundary condition is prescribed.

We assume that \eqref{kkt} is satisfied point-wisely as follows:
\begin{equation}\label{full:dis}
\begin{split}
&\partial_t u_h(\bz,t) +\mathcal{L}_hu_h+\mathcal{N}_h(u_h)=\lambda_h(\bz,t)g'(u_h),\quad\forall {\bm z} \in \Sigma_h,\\
&\lambda_h(\bz,t) \ge 0,\; g(u_h(\bz,t))\ge 0,\;\lambda_h(\bz,t) g(u_h(\bz,t))=0,\quad\forall {\bm z} \in \Sigma_h,
\end{split}
\end{equation}
with the Dirichlet boundary condition to be satisfied point-wisely if the original problem includes Dirichlet boundary condition at part or all of boundary.
  The above scheme includes  finite difference schemes, collocation schemes, or  Galerkin type spatial discretization with a Lagrangian basis.

Denote $\delta t$ the time step, and $t^n=n\delta t$ for $n=0,1,2,\cdots, \frac{T}{\delta t}$ where $T$ is the final computational time.  Our schemes consist of two steps: in the first step, we use a generic time discretization, which can be implicit, explicit or implicit-explicit, to find an intermediate solution $\tilde{u}_h^{n+1}$ which is usually not bound preserving; then we introduce a Lagrange multiplier  $\lambda_h^{n+1}(\bz)$ to determine a bound preserving ${u}_h^{n+1}$, which is a correction to $\tilde{u}_h^{n+1}$. We shall first construct bound preserving schemes which do not necessarily preserve mass, then we introduce a simple modification which allows us to construct  bound preserving schemes which can also preserve mass.

For the sake of clarity,  we shall restrict ourselves to constructed schemes based on the implicit-explicit (IMEX) type time discretization since they are most commonly used for parabolic type systems. It is straightforward to extend the approach below to schemes based on other types of time discretization.

\subsection{A class of multistep IMEX schemes}
We construct below $k$-th order bound-preserving  schemes for \eqref{full:dis} based on backward difference formula (BDF) for the time derivative and Adams-Bashforth extrapolation  by using a predictor-corrector approach. 

 In order to describe the scheme, we define a sequence $\{\alpha_k\}$, and with a slight abuse of notation.  For any function $v$, we use 
$A_k(v^n)$ and $B_{k-1}(v^n)$ to denote two operators depending on $(v^n, \cdots, v^{n-k+1})$ as follows:

\noindent {\bf $k=1$:}
\begin{equation}\label{eq:bdf1}
\alpha_1=1, \quad A_1(v^n)=v^n,\quad 
B_0(v^n)=0;
\end{equation}
\noindent {\bf $k=2$:}
\begin{equation}\label{eq:bdf2}
\alpha_2=\frac{3}{2}, \quad A_2(v^n)=2v^n-\frac{1}{2}v^{n-1},\quad B_1(v^n)=v^n;
\end{equation}
\noindent {\bf $k=3$:}
\begin{equation}\label{eq:bdf3}
\begin{split}
\alpha_3=\frac{11}{6}, &\quad A_3(v^n)=3v^{n}-\frac{3}{2}v^{n-1}+\frac{1}{3}v^{n-2},\\& B_2(v^n)=2v^n-v^{n-1}.
\end{split}
\end{equation}
The formula for $k=4,5,6$ can be derived similarly with Taylor expansions. 

We assume that $u_h^j, j=0,1,\cdots,k-1$ are properly initialized. Then

\noindent{\bf Step 1:} (Predictor) solve $\tilde{u}_h^{n+1}$ from
\begin{equation}\label{high:bound:lag:1}
\begin{split}
&\frac{\alpha_k\tilde{u}_h^{n+1}(\bz)-A_k(u_h^n(\bz))}{\delta t} + \mathcal{L}_h\tilde u_h^{n+1}(\bz)+\mathcal{N}_h(B_{k}(u_h^n(\bz)))=B_{k-1}(\lambda_h^ng'(u_h^n(\bz))), \quad\forall {\bm z} \in \Sigma_h;
\end{split}
\end{equation}
 {\bf Step 2:} (Corrector) solve $u_h^{n+1}$ and $\lambda_h^{n+1}$ from
 \begin{subequations}\label{high:bound:lag:2}
\begin{eqnarray}
&&\frac{\alpha_k(u_h^{n+1}(\bz)-\tilde{u}_h^{n+1}(\bz))}{\delta t} =\lambda_h^{n+1}(\bz)g'(u_h^{n+1}(\bz))-B_{k-1}(\lambda_h^n(\bz)g'(u_h^n(\bz))),\label{high:bound:lag:2a}\\
&&g(u_h^{n+1}(\bz)) \ge 0,\;\lambda_h^{n+1}(\bz)\ge 0,\; \lambda_h^{n+1}(\bz)g(u_h^{n+1}(\bz))=0,\quad\forall {\bm z} \in \Sigma_h.\label{high:bound:lag:2b}
\end{eqnarray}
\end{subequations}
 The second step     can be solved point-wisely as follows. 
 We  denote
\begin{equation}\label{defEta1}
\eta_h^{n+1}:=-\frac{\delta t}{\alpha_k}B_{k-1}(\lambda_h^ng'(u_h^n)),
\end{equation} and
 rewrite  \eqref{high:bound:lag:2a} as 
\begin{equation*}\label{por:lag:2d}
\frac{\alpha_k( u_h^{n+1}(\bz)-(\tilde{u}_h^{n+1}(\bz)+\eta^{n+1}_h(\bz)))}{\delta t} =\lambda_h^{n+1}(\bz)g'(u_h^{n+1}(\bz)). 
\end{equation*}
We find from the above and \eqref{high:bound:lag:2b} that
\begin{equation}
\begin{split}
 &(u_h^{n+1}(\bz),\lambda_h^{n+1}(\bz))=\left\{
\begin{array}{rcl}
(\tilde{u}_h^{n+1}(\bz)+\eta^{n+1}_h(\bz),0)       &      & {\mbox{if} \quad a      <  \tilde{u}_h^{n+1}(\bz)} +\eta^{n+1}_h(\bz) < b\\
(a,\frac{a-(\tilde{u}_h^{n+1}(\bz)+\eta^{n+1}_h(\bz))}{\frac{\delta t}{\alpha_k}g'(a)} )    &      &  {\mbox{if} \quad \tilde{u}_h^{n+1}(\bz)+\eta^{n+1}_h(\bz) \leq a}\\
(b,\frac{b-(\tilde{u}_h^{n+1}(\bz)+\eta^{n+1}_h(\bz))}{\frac{\delta t}{\alpha_k}g'(b)})    &      &  {\mbox{if} \quad \tilde{u}_h^{n+1}(\bz)  }+\eta^{n+1}_h(\bz) \ge b
\end{array} \right.,\;\forall \bz\in \Sigma_h.
\end{split}
\end{equation}

\begin{remark}\label{rem:cut:2}

It is obvious that the above scheme is a $k$-th order approximation to  \eqref{full:dis}. We would like to point out that  it is also a $k$-th order (in time) approximation plus the spatial discretization error  to  \eqref{strong}.

On the other hand, if we replace $B_{k-1}(\lambda_h^ng'(u_h^n))$ in the above scheme by zero, then it is easy to see that the second step is equivalent to the  simple cut-off approach, which is a first-order approximation to  \eqref{full:dis}.  However, it is easy to see that the error in maximum norm by the cut-off approach is smaller than the error by the corresponding semi-implicit scheme, therefore, the cut-off approach  is also  a $k$-th order (in time) approximation  plus the spatial discretization error to \eqref{strong}. 

\end{remark}

\subsection{Mass conservation}

A drawback of the schemes  \eqref{high:bound:lag:1}- \eqref{high:bound:lag:2} is that they do not preserve mass if the exact solution does. 

We present below a simple modification which enables mass conservation. More precisely, we introduce another Lagrange multiplier $\xi^{n+1}_h$, which is independent of spatial variables, to enforce the mass conservation in the second step. 

The first step is still exactly the same as \eqref{high:bound:lag:1}.

\noindent{\bf Step 1} (predictor): solve $\tilde u_h^{n+1}$ from
\begin{equation}
	\begin{split}
	\frac{\alpha_k \tilde{u}_h^{n+1}(\bz)-A_k(u_h^n(\bz))}{ \delta t} &+ \mathcal{L}_h \tilde u_h^{n+1}(\bz)+ \mathcal{N}_h (B_{k}(u_h^n(\bz)))\\
	&=B_{k-1}(\lambda_h^n(\bz)g'(u_h^n(\bz)))+B_{k-1}(\xi_h^n),\;\;\forall \bz \in \Sigma_h.
\end{split}\label{mass:por:lag:1}
\end{equation}

We introduce another Lagrange multiplier $\xi^{n+1}_h$ in the second step to enforce the mass conservation.

{\bf Step 2} (corrector): solve $(u_h^{n+1},\lambda_h^{n+1})$ from
\begin{subequations}\label{mass:por:lag:2}
	\begin{align}
		&\frac{\alpha_k( u_h^{n+1}(\bz)-\tilde{u}_h^{n+1}(\bz))}{\delta t} =\lambda_h^{n+1}(\bz)g'(u_h^{n+1}(\bz))\label{mass:por:lag:2a}\\
		&\hskip 3cm -B_{k-1}(\lambda_h^n(\bz)g'(u_h^n(\bz)))+\xi^{n+1}_h-B_{k-1}(\xi_h^n),\quad\forall {\bm z} \in \Sigma_h,\nonumber\\
		&\lambda_h^{n+1}(\bz)\ge 0,\; g(u_h^{n+1}(\bz))\ge 0,\; \lambda_h^{n+1}(\bz)g(u_h^{n+1}(\bz))=0,\;\quad\forall {\bm z} \in \Sigma_h,\label{mass:por:lag:2b}\\
		&( u_h^{n+1},1)_h=(u_h^{n},1)_h,\label{mass:por:lag:2c}
	\end{align}
\end{subequations}
where $(\cdot,\cdot)_h$ is a discrete inner product.

In order to solve the above system, we  denote
\begin{equation}\label{defEta}
	\eta_h^{n+1}:=\frac{\delta t}{\alpha_k}(\xi^{n+1}_h-B_{k-1}(\xi_h^n)-B_{k-1}(\lambda_h^ng'(u_h^n))),
\end{equation} and
rewrite  \eqref{mass:por:lag:2a} as 
\begin{equation}\label{mass:por:lag:2d}
	\frac{\alpha_k( u_h^{n+1}(\bz)-(\tilde{u}_h^{n+1}(\bz)+\eta^{n+1}_h(\bz)))}{\delta t} =\lambda_h^{n+1}(\bz)g'(u_h^{n+1}(\bz)). 
\end{equation}
Hence, assuming $\xi^{n+1}_h$ is known,  we find from the above and \eqref{mass:por:lag:2b}  that
\begin{equation}\label{mass:bound:sol}
	\begin{split}
		&(u_h^{n+1}(\bz),\lambda_h^{n+1}(\bz))=\left\{
		\begin{array}{rcl}
			(\tilde{u}_h^{n+1}(\bz)+\eta^{n+1}_h(\bz),0)       &      & {\mbox{if} \quad a      <  \tilde{u}_h^{n+1}(\bz)} +\eta^{n+1}_h(\bz) < b\\
			(a,\frac{a-(\tilde{u}_h^{n+1}(\bz)+\eta^{n+1}_h(\bz))}{\frac{\delta t}{\alpha_k}g'(a)} )    &      &  {\mbox{if} \quad \tilde{u}_h^{n+1}(\bz)+\eta^{n+1}_h(\bz) \leq a}\\
			(b,\frac{b-(\tilde{u}_h^{n+1}(\bz)+\eta^{n+1}_h(\bz))}{\frac{\delta t}{\alpha_k}g'(b)})    &      &  {\mbox{if} \quad \tilde{u}_h^{n+1}(\bz)  }+\eta^{n+1}_h(\bz) \ge b
		\end{array} \right.,\;\forall \bz\in \Sigma_h.
	\end{split}
\end{equation}

It remains  to determine $\xi^{n+1}_h$. 

Denote
\begin{equation}\label{sigmaabc}
\begin{split}
 &{}^a\Sigma_h(\xi)=\{z\in \Sigma_h:\,\tilde{u}_h^{n+1}(\bz)+\delta t\xi\le a\},\\
 &{}^a\Sigma_h^b(\xi)=\{z\in \Sigma_h:\,a<\tilde{u}_h^{n+1}(\bz)+\delta t\xi< b\},\\
& \Sigma_h^b(\xi)=\{z\in \Sigma_h:\,\tilde{u}_h^{n+1}(\bz)+\delta t\xi\ge b\}.
 \end{split}
 \end{equation}
 Then,  thanks to \eqref{mass:bound:sol}, the discrete mass conservation \eqref{mass:por:lag:2c}
 can be rewritten as
\begin{equation}\label{massn}
\sum_{z\in {}^a\Sigma_h^b(\eta^{n+1}_h)} (\tilde{u}_h^{n+1}(\bz)+\delta t\eta^{n+1}_h)\omega_z+\sum_{z\in \Sigma_h^b(\eta^{n+1}_h)} b\,\omega_z +\sum_{z\in {}^a\Sigma_h(\eta^{n+1}_h)}a\,\omega_z=(u_h^n,1)_h.
\end{equation}
Setting
\begin{equation}\label{defF}
	\begin{split}
	G_n(\eta)&:=\sum_{z\in {}^a\Sigma_h^b(\eta)} (\tilde{u}_h^{n+1}(\bz)+\delta t\eta)\omega_z+\sum_{z\in \Sigma_h^b(\eta)} b\,\omega_z +\sum_{z\in {}^a\Sigma_h(\eta)}a\,\omega_z-(u_h^n,1)_h,\\
	F_n(\xi)&:=G_n\big( \frac{\delta t}{\alpha_k}(\xi-B_{k-1}(\xi_h^n)-B_{k-1}(\lambda_h^ng'(u_h^n))) \big),
	\end{split}
\end{equation}
we find from  the above and \eqref{massn} that  $\xi^{n+1}_h$ is a solution to the nonlinear algebraic equation 
$F_n(\xi)=0.$
Since $F_n'(\xi)$ may not exist and is difficult to compute if it exists, instead of the Newton iteration, we can use the following secant method:
\begin{equation}\label{newton}
\xi_{k+1}=\xi_k -\frac{F_n(\xi_k)(\xi_k-\xi_{k-1})}{F_n(\xi_k)-F_n(\xi_{k-1})}.
\end{equation}
Since $\xi^{n+1}_h$ is an approximation to zero, we can choose  $\xi_0=0$ and $\xi_1=O(\delta t)$. In all our experiments,  \eqref{newton} converges in a few iterations so that the cost is negligible.

 Once $\xi^{n+1}_h$ is known, we can update $(u_h^{n+1},\lambda^{n+1}_h)$ with \eqref{mass:bound:sol}.

\begin{remark}
It is usually very difficult to construct mass conserved IMEX schemes using the simple cut-off approach. However, 
	replacing $B_{k-1}(\lambda_h^n(\bz)g'(u_h^n(\bz)))$ in \eqref{mass:por:lag:1}-\eqref{mass:por:lag:2} by zero, we obtain a mass conserved $k$th-order IMEX cut-off scheme. This is one of the advantages of reformulating the cut-off approach with the operator splitting approach.
\end{remark}

\section{Stability results}
While the schemes constructed in the last section automatically ensure the $L^\infty$ bound for $\{u_h^n\}$, it does not imply any bound on the energy norm $<\mathcal{L}\cdot,\cdot>$. In this section, we shall use the energy estimates to derive a bound on the energy norm for $\{\tilde u_h^n\}$ as well as a bound on the Lagrange multiplier. 

To fix the idea, we assume that $ \mathcal{L}$ is a second-order unbounded positive self-adjoint operator in $L^2(\Omega)$ with domain $D( \mathcal{L})$, and that the nonlinear term can be written as follows:
\begin{equation}\label{assump1}
\begin{split}
&  \mathcal{N}(u)=f_1(u)+\nabla\cdot f_2(u), \quad\text{with }\;f_1(0)=f_2(0)=0, \\
& \text{ and }\, f_1,\, f_2 \;\text{   are locally Lipchitz semi-linear functions}.
\end{split}
\end{equation}
Without loss of generality, we assume that $ab \le 0$. Otherwise, we can always find a constant $C$ such that $(a+C)(b+C) \le 0$ and consider the equation for $v=u+C$.
Since $ab \le 0$, we have $0\in (a,b)$. Hence,
\eqref{assump1} implies in particular
\begin{equation}\label{assump0}
 |f_1(u)|=|f_1(u)-f_1(0)|\le C_1 |u|,\;  |f_2(u)|=|f_2(u)-f_2(0)|\le C_2 |u|\quad \text{ if }\; a\le u\le b.
\end{equation}
We observe that the nonlinearities  in common nonlinear parabolic equations do satisfy \eqref{assump1}, see in particular some specific examples given in Section 5.

We shall also interpret the first step of the schemes, \eqref{high:bound:lag:1} and \eqref{mass:por:lag:1}, in  a Galerkin formulation. More precisely, let $X_h\subset X$ be a subspace  with Lagrangian basis functions  on $\Sigma_h$. We 
define a discrete inner product on  $\Sigma_h=\{\bm{z}\}$ in $\bar \Omega$: 
\begin{equation}\label{numint}
	(u,v)_h=\sum_{\bm{z}\in \Sigma_h} \beta_{\bm{z}}u(\bm{z})v(\bm{z}), 
\end{equation}
where we require that the weights $\beta_{\bm{z}}>0$. We also denote the induced norm by $\|u\|=(u,u)_h^{\frac 12}$, and we assume that this norm is equivalent to the $L^2$ norm for functions in $X_h$. 
We denote by $< \mathcal{L}_h u_h,v_h>$ the bilinear form on $X_h\times X_h$ based on  the discrete inner product after suitable integration by part, and
we  assume that 
\begin{equation}\label{assump2}
C_0\|\nabla u_h\|^2\le <\mathcal{L}_h u_h,u_h> \quad\forall u_h\in X_h,
\end{equation}
with $C_0>0$, which is satisfied by many common spatial discretizations. Hereafter, we shall use  $C$ and $C_i$ to denote  generic positive constants which are independent of $\delta t$ and $h$.

We shall only consider a second-order scheme with mass conservation  in this section. It is clear that similar bounds can be derived  for second-order scheme without mass conservation, and for  the first-order schemes, but bounds for higher-order schemes are still elusive.  For clarity, we rewrite the second-order version of  \eqref{mass:por:lag:1}- \eqref{mass:por:lag:2} as:

\noindent{\bf Step 1} (predictor): Find  $\tilde u_h^{n+1}\in X_h$ such that, for  $\forall v_h\in X_h$
\begin{align}\label{2nd:por:lag:1}
	(\frac{3 \tilde{u}_h^{n+1}-4u_h^n+u_h^{n-1}}{ 2\delta t},v_h)_h+ <\mathcal{L}_h \tilde u_h^{n+1},v_h>+
(f_1(u_h),v_h)_h-(f_2(u_h),\nabla v_h)_h=(\lambda_h^ng'(u_h^n)+\xi_h^n,v_h)	_h;
\end{align}

\noindent{\bf Step 2} (corrector): Find $u_h^{n+1},\,\lambda_h^{n+1},\,\xi^{n+1}_h$ from
\begin{subequations}\label{2nd:por:lag:2}
	\begin{align}
		&\frac{3( u_h^{n+1}(\bz)-\tilde{u}_h^{n+1}(\bz))}{2\delta t} =\lambda_h^{n+1}(\bz)g'(u_h^{n+1}(\bz)) -\lambda_h^n(\bz)g'(u_h^n(\bz))+\xi^{n+1}_h-\xi_h^n,\quad\forall {\bm z} \in \Sigma_h,\label{2nd:por:lag:2a}\\
		&\lambda_h^{n+1}(\bz)\ge 0,\; g(u_h^{n+1}(\bz))\ge 0,\; \lambda_h^{n+1}(\bz)g(u_h^{n+1}(\bz))=0,\;\quad\forall {\bm z} \in \Sigma_h,\label{2nd:por:lag:2b}\\
		&( u_h^{n+1},1)_h=(u_h^{n},1)_h;\label{2nd:por:lag:2c}
	\end{align}
\end{subequations}
and we assume that $\tilde{u}_h^0$ and ${u}_h^0$ are computed with the first-order scheme \eqref{mass:por:lag:1}-\eqref{mass:por:lag:2} with $k=1$.

\begin{theorem}
We assume \eqref{assump1},  \eqref{assump0} and  \eqref{assump2}.
	Then, for  the scheme \eqref{2nd:por:lag:1}-\eqref{2nd:por:lag:2},  if the generic scheme in \eqref{mass:por:lag:1} is mass conservative, i.e., 
	\begin{equation}\label{massassump4}
		 <\mathcal{L}_h \tilde u_h^{n+1},1>+
		(f_1(u_h),1)_h-(f_2(u_h),\nabla 1)_h=0,
	\end{equation}
then, 	we have 
	\begin{equation*}\label{mass:th:1}
		\begin{split}
			4\|u_h^{m}\|^2
			+ \|2 u_h^{m}- u_h^{m-1}\|^2&+\frac 43\delta t^2\|\lambda_h^{n+1}g'(u_h^{m})+\xi_h^{m}\|^2 \\&+2\delta t\sum\limits_{n=0}^{m-1}C_0\|\nabla  \tilde{u}_h^{n+1}\|^2
			\leq C(T)\| u_h^0\|^2,\;\forall 1\le m\le T/{\delta t}.
		\end{split}
	\end{equation*}
	
\end{theorem}
\begin{proof}
		Choosing $v_h=4\delta t\tilde{u}_h^{n+1}$ in \eqref{2nd:por:lag:1}, using  the assumption \eqref{assump2}, we obtain
	\begin{equation}\label{bdf2:stab:1}
		\begin{split}
			&(3\tilde{u}_h^{n+1}-4 u_h^n+ u_h^{n-1},2\tilde{ u}_h^{n+1})_h
			+4\delta tC_0\|\nabla  \tilde{u}_h^{n+1}\|^2\\
			&+4\delta t(f_1(2u_h^{n}-u_h^{n-1}),\tilde{u}_h^{n+1})_h-4\delta t(f_2(2u_h^{n}-u_h^{n-1}),\nabla  \tilde{u}_h^{n+1})_h
			\leq 4\delta t(\lambda_h^n g'(u_h^n)+\xi_h^n,\tilde{u}_h^{n+1})_h.
		\end{split}
	\end{equation}
	We start by dealing with   the first term in \eqref{bdf2:stab:1}.
	\begin{equation}\label{bdf2:stab:2}
		\begin{split}
			&(3\tilde{u}_h^{n+1}-4 u_h^n+ u_h^{n-1},2\tilde{u}_h^{n+1})_h=
			2(3 u_h^{n+1}-4u_h^n+ u_h^{n-1}, u_h^{n+1})_h \\& +6(\tilde{u}_h^{n+1}- u_h^{n+1},\tilde{u}_h^{n+1})_h+ 2(3 u^{n+1}_h-4 u_h^n+ u_h^{n-1},\tilde{u}_h^{n+1}- u_h^{n+1})_h.
		\end{split}
	\end{equation}
	For the terms on the righthand side of  \eqref{bdf2:stab:2}, we have
	\begin{equation}\label{bdf2:stab:3}
		\begin{split}
			&2(3 u_h^{n+1}-4 u_h^n+ u_h^{n-1}, u_h^{n+1})_h=
			\| u_h^{n+1}\|^2-\| u_h^n\|^2\\&+\|2 u_h^{n+1}- u_h^n\|^2-\|2 u_h^n- u_h^{n-1}\|^2+\| u_h^{n+1}-2 u_h^n+ u_h^{n-1}\|^2;
		\end{split}
	\end{equation}
	\begin{equation}\label{bdf2:stab:5}
		6(\tilde{u}_h^{n+1}- u_h^{n+1},\tilde{u}_h^{n+1})_h=3(\|\tilde{u}_h^{n+1}\|^2-\| u_h^{n+1}\|^2+\|\tilde{u}_h^{n+1}- u_h^{n+1}\|^2);
	\end{equation}
	and 
	\begin{equation}\label{bdf2:stab:4}
		\begin{split}
			2(3 u_h^{n+1}&-4 u_h^n+ u_h^{n-1},\tilde{u}_h^{n+1}- u_h^{n+1})_h\\
			&=2(u_h^{n+1}-2 u_h^n+ u_h^{n-1},\tilde{u}_h^{n+1}- u_h^{n+1})_h+
			4( u_h^{n+1}- u_h^n,\tilde{u}_h^{n+1}- u_h^{n+1})_h\\
			&\ge -\|u_h^{n+1}-2 u_h^n+ u_h^{n-1}\|^2 -\|\tilde{u}_h^{n+1}- u_h^{n+1}\|^2+4( u_h^{n+1}- u_h^n,\tilde{u}_h^{n+1}- u_h^{n+1})_h.
		\end{split}
	\end{equation}
	The last term in the above needs a special treatment. Using \eqref{2nd:por:lag:2a} and the fact that
	$(u_h^{n+1}- u_h^n,1)_h=0$, we can write
	\begin{equation}\label{larger}
		\begin{split}
			&4( u_h^{n+1}- u_h^n,\tilde{ u}_h^{n+1}- u_h^{n+1})_h=-\frac{8\delta t}{3}(u_h^{n+1}- u_h^n,\lambda_h^{n+1}g'(u_h^{n+1})-\lambda_h^ng'(u_h^n)+\xi_h^{n+1}-\xi_h^n)_h\\&
			=-\frac{8\delta t}{3}(u_h^{n+1}- u_h^n,\lambda_h^{n+1}g'(u_h^{n+1})-\lambda_h^ng'(u_h^n))_h -\frac{8\delta t}{3}(\xi_h^{n+1}-\xi_h^n)(u_h^{n+1}- u_h^n,1)_h
			\\&=-\frac{8\delta t}{3}( u_h^{n+1}- u_h^n,\lambda_h^{n+1}g'(u_h^{n+1}))_h -\frac{8\delta t}{3}( u_h^{n}- u_h^{n+1},\lambda_h^{n}g'(u_h^{n}))_h:=I_1 +I_2.
		\end{split}
	\end{equation}
Thanks to $\lambda_h^{n+1}(\bz)g(u_h^{n+1}(\bz))=0$,
	we obtain
	\begin{equation*}\label{stab:i1}
		\begin{split}
			I_1&= -\frac{8\delta t}{3}(\lambda_h^{n+1}, (u_h^{n+1}- u_h^n)(a+b-2u_h^{n+1})-g(u_h^{n+1}))_h 
			\\&= -\frac{8\delta t}{3} (\lambda_h^{n+1}, -(u_h^{n+1})^2+ab+2u_h^nu_h^{n+1}-(a+b)u_h^n )_h
			\\&=\frac{8\delta t}{3} ( \lambda_h^{n+1}, (u_h^{n+1}-u_h^n)^2)_h-\frac{8\delta t}{3} ( \lambda^{n+1}_h, (u_h^n-a)(u_h^n-b))_h\ge 0,
		\end{split}
	\end{equation*}
	where we used the facts that $a\leq u_h^n \leq b$ and $\lambda_h^{n+1}\ge 0$.
	Similarly, we use $\lambda_h^{n}(\bz)g(u_h^{n}(\bz))=0$  to derive
	\begin{equation*}\label{stab:i2}
		\begin{split}
			I_2&= -\frac{8\delta t}{3}(\lambda_h^{n}, (u_h^{n}- u_h^{n+1})g'(u_h^{n})-g(u_h^{n}))_h
			\\&= -\frac{8\delta t}{3} (\lambda_h^{n},-(u_h^{n})^2+ab+2u_h^nu_h^{n+1}-(a+b)u_h^{n+1} )_h
			\\&=\frac{8\delta t}{3}( \lambda_h^{n},(u_h^{n+1}-u_h^n)^2)_h-(\lambda_h^n, (u_h^{n+1}-a)(u_h^{n+1}-b))_h\ge 0,
		\end{split}
	\end{equation*}
	where we used again the facts that  $\lambda^n_h \ge 0$ and $a \leq u_h^{n+1}\leq b$.
	We derive from the last two inequalities that
	\begin{equation}\label{mass:larger}
		\begin{split}
			&4( u_h^{n+1}- u_h^n,\tilde{ u}_h^{n+1}- u_h^{n+1})_h=-\frac{8\delta t}{3}( u_h^{n+1}- u_h^n,\lambda_h^{n+1}g'(u_h^{n+1})-\lambda_h^ng'(u_h^n))_h\ge 0.
		\end{split}
	\end{equation}
	Combining the above inequalities in \eqref{bdf2:stab:2}, we find
	\begin{equation}\label{bdf2:stab:2b}
		\begin{split}
			(3\tilde{u}_h^{n+1}-4 u_h^n+ u_h^{n-1},2\tilde{u}_h^{n+1})_h&\ge  \| u_h^{n+1}\|^2-\| u_h^n\|^2+\|2 u_h^{n+1}- u_h^n\|^2-\|2 u_h^n- u_h^{n-1}\|^2 \\
			&+3( \| \tilde u_h^{n+1}\|^2-\| u_h^{n+1}\|^2) + 2\|\tilde u_h^{n+1}- u_h^{n+1}\|^2.
		\end{split}
	\end{equation}
	
	Next, we rewrite \eqref{2nd:por:lag:2a} as
	\begin{equation}\label{eq:proj}
		3 u_h^{n+1}-2\delta t(\lambda_h^{n+1}g'(u_h^{n+1})+\xi_h^{n+1})=3\tilde{u}_h^{n+1}-2\delta t(\lambda_h^ng'(u_h^n)+\xi_h^n).
	\end{equation}
	Taking the discrete inner product of each side of the equation \eqref{eq:proj} with itself, dividing by $3$, we obtain
	\begin{equation}\label{mass:bdf2:stab:2} 
		\begin{split}
			3\| u_h^{n+1}\|^2&-4\delta t(u_h^{n+1},\lambda_h^{n+1}g'(u_h^{n+1})+\xi_h^{n+1})_h+\frac 43\delta t^2\|\lambda_h^{n+1}g'(u_h^{n+1})+\xi_h^{n+1}\|^2 \\&= 3\|\tilde{u}_h^{n+1}\|^2-4\delta t(\tilde{u}_h^{n+1},\lambda_h^ng'(u_h^n)+\xi_h^n)_h + \frac 43\delta t^2\|\lambda_h^ng'(u_h^n)+\xi_h^n\|^2.
		\end{split}
	\end{equation}
Note that we can interpret  \eqref{2nd:por:lag:1}  pointwisely as
\begin{align}\label{2nd:por:lag:1b}
\frac{3 \tilde{u}_h^{n+1}(\bz)-4u_h^n(\bz)+u_h^{n-1}(\bz)}{ 2\delta t}+ \mathcal{L}_h \tilde u_h^{n+1}(\bz)+
\mathcal{N}_h (2u_h^n(\bz)-u_h^{n-1}(\bz))=\lambda_h^n(\bz)g'(u_h^n(\bz))+\xi_h^n,\;\forall \bz\in \Sigma_h,
\end{align}
where $ \mathcal{N}_h$ is defined by $(\mathcal{N}_h (u_h),v_h)_h=(f_1(u_h),v_h)_h-(f_2(u_h),\nabla v_h)_h$. 
	Summing  up \eqref{2nd:por:lag:1b} and \eqref{2nd:por:lag:2a}, we obtain
	\begin{align}\label{mass:prod}
		\frac{3 {u}_h^{n+1}(\bz)-4u_h^n(\bz)+{u}_h^{n-1}(\bz)}{2 \delta t}+ \mathcal{L}_h \tilde u_h^{n+1}(\bz)+ \mathcal{N}_h (2 u_h^n(\bz)- u_h^{n-1}(\bz))=\lambda_h^{n+1}(\bz)g'(u_h^{n+1}(\bz))+\xi_h^{n+1},\;\forall \bz\in \Sigma_h.
	\end{align}
	Taking the discrete inner product of \eqref{mass:prod} with $1$ on both sides, using \eqref{2nd:por:lag:2c} and \eqref{massassump4}, we obtain
	\begin{equation}
		(\lambda_h^{n+1}g'(u_h^{n+1})+\xi_h^{n+1}, 1)_h=0,
	\end{equation}
	which implies that 
	\begin{equation}\label{mass:stab:5}
		\xi_h^{n+1}=-\frac{(\lambda_h^{n+1}g'(u_h^{n+1}),1)_h}{|\Omega|}
		=-\frac{(\lambda_h^{n+1},a+b-2u_h^{n+1})_h}{|\Omega|},
	\end{equation}
	where $|\Omega|:=(1,1)_h=\Sigma_{\bz\in \Sigma_k}\beta_{\bz}>0$.
	
	It remains to show that the second term  of 
	\eqref{mass:bdf2:stab:2} is non negative. Using  the fact that $\lambda_h^{n+1}(\bz)g(u_h^{n+1}(\bz))=0$, we have
	\begin{equation*}\label{mass:stab:6}
		\begin{split}
			-4\delta t(u_h^{n+1},&\lambda_h^{n+1}g'(u_h^{n+1})+\xi_h^{n+1})_h
			=-4\delta t(\lambda_h^{n+1},u_h^{n+1}g'(u_h^{n+1})-g(u_h^{n+1}))_h-4\delta t\xi_h^{n+1}(u_h^{n+1},1)_h
			\\&=-4\delta t(\lambda_h^{n+1}, ab-(u_h^{n+1})^2)_h+\frac{4\delta t}{|\Omega|}(\lambda_h^{n+1},a+b-2u_h^{n+1})_h(u_h^{n+1},1)_h
			\\&=-4\delta t(\lambda_h^{n+1},-(u_h^{n+1}-\frac{(u_h^{n+1},1)_h}{|\Omega|})^2+(\frac{(u_h^{n+1},1)_h}{|\Omega|}-a)(\frac{(u_h^{n+1},1)_h}{|\Omega|}-b))_h.
		\end{split}
	\end{equation*}
	Since $a\le u_h^{n+1} \le b$, we have
	\begin{equation}
		(\frac{(u_h^{n+1},1)_h}{|\Omega|}-a)(\frac{(u_h^{n+1},1)_h}{|\Omega|}-b) \le 0,
	\end{equation}
	which, together with  $\lambda_h^{n+1} \ge 0$, implies that
	\begin{equation}\label{mass:leq}
		-4\delta t(u_h^{n+1},\lambda_h^{n+1}g'(u_h^{n+1})+\xi_h^{n+1})_h \ge 0.
	\end{equation}
	Then, summing up \eqref{bdf2:stab:1} with \eqref{mass:bdf2:stab:2}, and using \eqref{mass:larger}, \eqref{bdf2:stab:2b} and\eqref{mass:leq},  after dropping some unnecessary terms, we obtain
	\begin{equation}\label{mass:stab:final:new}
		\begin{split}
			&4\| u_h^{n+1}\|^2-4\| u_h^n\|^2+\|2 u_h^{n+1}- u_h^n\|^2-\|2 u_h^n- u_h^{n-1}\|^2+ 2 \|\tilde u_h^{n+1}- u_h^{n+1}\|^2
			\\&+\frac 43\delta t^2(\|\lambda_h^{n+1}g'(u_h^{n+1})+\xi_h^{n+1}\|^2-\|\lambda_h^ng'(u_h^n)+\xi_h^n\|^2)
			+ 4\delta tC_0\|\nabla  \tilde{u}_h^{n+1}\|^2\\&\le -4\delta t(f_1(2u_h^{n}-u_h^{n-1}),\tilde{u}_h^{n+1})_h+4\delta t(f_2(2u_h^{n}-u_h^{n-1}),\nabla  \tilde{u}_h^{n+1})_h.
		\end{split}
	\end{equation}
	Using \eqref{assump2}, the two terms on the righthand side above can be bounded as follows:
	\begin{equation}\label{young:1}
		\begin{split}
			4\delta t(f_1(2u_h^{n}-u_h^{n-1}),\tilde{u}_h^{n+1})_h&=4\delta t(f_1(2u_h^{n}-u_h^{n-1}),\tilde{u}_h^{n+1}-u_h^{n+1})_h+ 4\delta t(f_1(2u_h^{n}-u_h^{n-1}),u_h^{n+1})_h\\
			& \leq 2\|\tilde{u}_h^{n+1}-u_h^{n+1}\|^2 +2C_1^2\delta t^2\|2u_h^{n}-u_h^{n-1}\|^2\\
			& + 2\delta t(C_1^2\|2u_h^{n}-u_h^{n-1}\|^2+\|u_h^{n+1}\|^2).
		\end{split}
	\end{equation}
	Similarly, we have
	\begin{equation}\label{young:2}
		\begin{split}
			4\delta t(f_2(2u_h^{n}-u_h^{n-1}),\nabla  \tilde{u}_h^{n+1})_h &\leq 2\delta t C_0\|\nabla  \tilde{u}_h^{n+1}\|^2 +\frac{2\delta t}{C_0}\|f_2(2u_h^{n}-u_h^{n-1})\|^2
			\\&\leq  2\delta t C_0\|\nabla  \tilde{u}_h^{n+1}\|^2 +\frac{2\delta tC_2^2}{C_0}\|2u_h^{n}-u_h^{n-1}\|^2.
		\end{split}
	\end{equation}
	Combining \eqref{mass:stab:final:new}, \eqref{young:1} and \eqref{young:2}, we obtain
	\begin{equation}\label{mass:sum}
		\begin{split}
			4\| u_h^{n+1}\|^2-4\| u_h^n\|^2&+\|2 u_h^{n+1}- u_h^n\|^2-\|2 u_h^n- u_h^{n-1}\|^2
			\\&+\frac 43\delta t^2(\|\lambda_h^{n+1}g'(u_h^{n+1})+\xi_h^{n+1}\|^2-\|\lambda_h^ng'(u_h^n)+\xi_h^n\|^2)
			+ 2\delta tC_0\|\nabla  \tilde{u}_h^{n+1}\|^2\\&\le C\delta t\|2u_h^{n}-u_h^{n-1}\|^2 +2\delta t\|u_h^{n+1}\|^2,\quad\forall n\ge 1.
		\end{split}
	\end{equation}
	For $n=0$, we use a first-order scheme, namely  \eqref{mass:por:lag:1}-\eqref{mass:por:lag:2} with $k=1$, to compute $\tilde u_h^1$ and $u_h^1$.  Using a similar  (but much simplified) procedure as above, we can obtain 
	\begin{equation}\label{mass:sum0}
		\begin{split}
			\| u_h^{1}\|^2-\| u_h^0\|^2 &
			+\delta t^2(\|\lambda_h^{1}g'(u_h^{1})+\xi_h^{1}\|^2-\|\lambda_h^0g'(u_h^0)+\xi_h^0\|^2)
			+ 2\delta tC_0\|\nabla  \tilde{u}_h^{1}\|^2\\&\le C\delta t\|2u_h^{0}\|^2 +2\delta t\|u_h^{1}\|^2.
		\end{split}
	\end{equation}
	Finally summing up \eqref{mass:sum0} with \eqref{mass:sum} from $n=1$ to $n=m-1$, we obtain 
	\begin{equation*}\label{mass:stab:gron:1}
		\begin{split}
			&4\| u_h^{m}\|^2+\|2 u_h^{m}- u_h^{m-1}\|^2+\frac 43\delta t^2\|\lambda_h^{m}g'(u_h^m)+\xi_h^m\|^2+
			2\delta t \sum\limits_{n=1}^{m-1} C_0\|\nabla  \tilde{u}_h^{n+1}\|^2
			\\&\le \|2 u_h^1- u_h^{0}\|^2+4\| u_h^0\|^2+C\delta t \sum\limits_{n=0}^{m-1}\{\|2u_h^{n}-u_h^{n-1}\|^2 +\|u_h^{n+1}\|^2\}.
		\end{split}
	\end{equation*}
	Applying the discrete Gronwall lemma,  and using \eqref{mass:sum0},  we arrive at the desired result.
\end{proof}

\section{Error estimate}
The error analysis for the second-order scheme \eqref{2nd:por:lag:1}-\eqref{2nd:por:lag:2} with a general spatial discretization is very tedious and may obscure its essential difficulty. Therefore, we shall carry out a 
complete error analysis for a second-order bound preserving scheme with a hybrid spectral discretization that we shall describe below. To further  simplify the presentation,   we assume $\mathcal{L}=-\Delta$ with Dirichlet boundary conditions on $\Omega=(-1,1)^d\;(d=1,2,3). $

   We now describe some preliminaries for our  hybrid spectral discretization.
Let $P_N$ be the space of  polynomials of degree less than or equals to  $N$ in each direction, we set
\begin{equation}
	X=H^1_0(\Omega),\quad X_N=\{ v\in P_N:  v|_{\partial \Omega}=0\}.
\end{equation}
We define the projection operator $\Pi_N: X\rightarrow X_N $  by
\begin{equation}\label{de:proj}
	(\nabla(v-\Pi_Nv) ,\nabla v_N)=0,\;\forall v\in X,\; v_N \in X_N,
\end{equation}
and recall that for any $r\ge 1$, we have \cite{canuto2012spectral}
\begin{equation}\label{assump3}
	\|v-\Pi_N v\|_{H^s} \lesssim   N^{s-r}\|v\|_{H^r}, \quad \forall v\in H^r(\Omega)\cap X,\, (s=0,\, 1),
\end{equation}
where $\|\cdot\|_{H^r}$ denote the usual norm in $H^r(\Omega)$.

Let $L_N$ be the Legendre polynomial of degree $N$, and $\{x_k\}_{0\le k\le N}$ be the roots of $(1-x^2)L_N'(x)$, i.e., the   Legendre-Gauss-Lobatto points.
We set $\Sigma_N=\{x_k\}_{1\le k\le N-1}$ and  $\bar \Sigma_N=\{x_k\}_{0\le k\le N}$  if $d=1$,  $\Sigma_N=\{(x_k,x_i)\}_{1\le k,i\le N-1}$  and $\bar \Sigma_N=\{(x_k,x_i)\}_{0\le k,i\le N}$ if $d=2$ and   $\Sigma_N=\{(x_k,x_i,x_j)\}_{1\le k,i,j\le N-1}$ and $\bar \Sigma_N=\{(x_k,x_i,x_j)\}_{0\le k,i,j\le N}$ if $d=3$.  
We define the interpolation operator $I_N: C(\Omega)\rightarrow P_N$ by $(I_{N}u)(\bz)=u(\bz)$ for all $\bz\in \bar\Sigma_N$. Then, we also have  \cite{canuto2012spectral}
\begin{equation}\label{assump4}
	\|v-I_N v\|_{H^s} \lesssim   N^{s-r}\|v\|_{H^r}, \quad \forall v\in H^r(\Omega)\cap X,\, (s=0,\, 1).
\end{equation}

Let $(\cdot,\cdot)_N$  be the discrete inner product based on the Gauss-Lobatto quadrature, then it is well known that \cite{shen2011spectral}
\begin{equation}\label{disinner}
	\begin{split}
		(u_N,v_N)_N&=(u_N,v_N) \quad\forall u_N\cdot v_N \in P_{2N-1},\\
		\|v_N\|^2& \le (v_N,v_N)_N\le (2+1/N)	\|v_N\|^2\quad  \forall v_N\in {P_N}.
	\end{split}
\end{equation}

We observe that the bound preserving is enforced at the second step, so the first-step in the bound preserving schemes can be replaced by any other $k$-th order scheme.  We shall consider a second-order modified Crank-Nicholson scheme which is easier to analyze.
More precisely, we consider the following modified Crank-Nicholson scheme \cite{gottlieb2012stability} with a hybrid spectral discretization: find $u_N^{n+1}\in X_N$ such that for all $n\ge 1$,
\begin{equation}\label{spectral:sche}
	\begin{split}
		(\frac{\tilde u_N^{n+1}(\bz)- u_N^n(\bz)}{\delta t} ,v_N)_N&+ (\nabla \frac{3\tilde{u}_N^{n+1}(\bz)+\tilde u_N^{n-1}(\bz)}{4},\nabla v_N)\\&+(\mathcal{N}(\frac 32u_N^n(\bz)-\frac 12u_N^{n-1}(\bz)),v_N)= 0
		,\quad\forall v_N \in X_N;
	\end{split}
\end{equation}		
and find $u_N^{n+1} $, $\lambda_N^{n+1}$ such that
\begin{equation}\label{spectral:sche2}
	\begin{split}
		\frac{ u_N^{n+1}(\bz)-\tilde{ u}_N^{n+1}(\bz)}{\delta t} &=\lambda_N^{n+1}(\bz)g'(u_N^{n+1}(\bz)),\quad\forall {\bm z} \in \Sigma_N,\\
		&\lambda_N^{n+1}(\bz)\ge 0,\; g(u_N^{n+1}(\bz))\ge 0,\; \lambda_N^{n+1}(\bz) g(u_N^{n+1}(\bz))=0,\quad\forall {\bm z} \in \Sigma_N.
	\end{split}
\end{equation}
For $n=0$, we replace $\mathcal{N}(\frac 32u_N^n(\bz)-\frac 12u_N^{n-1}(\bz))$ in \eqref{spectral:sche} by 
$\mathcal{N}(u_N^n(\bz))$.

To simplify the notation, we shall use $u(t)$ to denote $u(\bx,t)$.   We denote
\begin{equation}
	\bar e_N^{n+1} =u(t^{n+1})-\Pi_N u(t^{n+1}),\; \hat e_N^{n+1}=\Pi_N u(t^{n+1}) - u_N^{n+1},\;
	\tilde e_N^{n+1}=\Pi_N u(t^{n+1})-\tilde u_N^{n+1}.
\end{equation}
Then, we have
\begin{equation}
	u(t^{n+1}) - u_N^{n+1}=\bar e_N^{n+1}+\hat e_N^{n+1},\; u(t^{n+1})-\tilde{u}_N^{n+1}=\bar e_N^{n+1}+\tilde e_N^{n+1}.
\end{equation}
Let $t^k=k\delta t$, $t^{k+\frac12}=\frac 12(t^{k+1}+t^k)$ and $u^{n+\frac 12}=\frac{u^{n+1}+u^n}{2}$. We denote
\begin{equation}
	\begin{split}
		&K_N^{n+\frac 12}=\frac{\bar e_N^{n+1}-\bar e_N^n}{\delta t},\\
		&T_N^{n+\frac 12}= -\Delta (u(t^{n+\frac 12})- \frac{3u(t^{n+1}) +u(t^{n-1})}{4}),\\
		&R_N^{n+\frac 12}=\partial_t u(t^{n+\frac12})-\frac{u(t^{n+1})-u(t^n)}{2},\\
		& J_N^{n+\frac 12}=u(t^{n+\frac 12})- (\frac 32 u(t^n)-\frac 12u(t^{n-1})).
	\end{split}
\end{equation}

\begin{theorem}
	Let $\tilde{u}_N^{n+1}, \, u_N^{n+1},\,\lambda_N^{n+1}$ be the solution of \eqref{spectral:sche}-\eqref{spectral:sche2}. Given $T\ge 0$, for some $l\ge 1$, assuming \eqref{assump1}-\eqref{assump0}, and the exact solution of \eqref{strong} $u(\bx,t) \in C^2([0,T],H^2(\Omega)) \cap C^1([0,T],H^l(\Omega))\cap C^3([0,T],L^2(\Omega))$, then we have the following error estimate:
	\begin{equation*}
		\begin{split}
	&	\|u(t^{m})-u_N^{m}\|^2+ \frac{\delta t}{4}\|\nabla(u(t^m)-\tilde{u}_N^{m})\|^2+ \frac{\delta t}{4}\|\nabla(u(t^{m-1})-\tilde{u}_N^{m-1})\|^2+\delta t^2\sum\limits_{n=1}^{m-1}\|\lambda_N^{n+1}g'(u_N^{n+1})\|^2_N\\
		&+\delta t\sum\limits_{n=1}^{m-1}\|\nabla(u(t^{n+1})-\tilde{u}_N^{n+1}+u(t^{n-1})-\tilde u^{n-1}_N)\|^2
		 \le  C(\delta t^4 +N^{-2l}), \quad \forall 2\leq m \leq \frac{T}{\delta t}.
		\end{split}
	\end{equation*}
\end{theorem}
\begin{proof}
	We derive from \eqref{strong} and \eqref{de:proj} that
	\begin{equation}\label{error:ex}
		(\partial_t  u, v_N)_N +(\nabla \Pi_N u, \nabla v_N)+(\mathcal{N}(u),v_N) =	\eps(v_N),\quad\forall v_N\in X_N,
	\end{equation}
	where 
	\begin{equation}
		\eps(v_N)= (\partial_t  u, v_N)_N - (\partial_t  u, v_N).
	\end{equation}
We find from  \eqref{disinner}, the definition of $I_N$  and \eqref{assump3}-\eqref{assump4} that
\begin{equation}\label{errdisinner}
	\begin{split}
	|	\eps(v_N)|&=	|	(u_t, v_N)_N - (u_t, v_N)|=|(u_t-\Pi_{N-1}u_t,v_N)_N +(\Pi_{N-1}u_t-u_t,v_N)|\\
		&=|(I_Nu_t-\Pi_{N-1}u_t,v_N)_N +(\Pi_{N-1}u_t-u_t,v_N)|\\
		&\le (3 \|I_Nu_t-\Pi_{N-1}u_t\| +\|\Pi_{N-1}u_t-u_t\|)\|v_N\|  \\
		&	\le  (3\|I_Nu_t-u_t\|+4\|u_t-\Pi_{N-1}u_t\| )\|v_N\| \le C N^{-l}\|v_N\|, \quad\forall v_N\in P_N.
	\end{split}	
\end{equation}
	Subtracting equation \eqref{error:ex}  from scheme \eqref{spectral:sche}, we obtain 
	\begin{equation}\label{error:01}
		\begin{split}
			&(\frac{\tilde{e}_N^{n+1}-\hat e_N^n}{\delta t},v_N ) _N
			+(\nabla \frac{3\tilde e_N^{n+1}+\tilde e_N^{n-1}}{4},\nabla v_N)+ (\mathcal{N}(u(t^{n+\frac 12}))-\mathcal{N}(\frac 32u_N^n-\frac 12 u_N^{n-1}),v_N) \\&=( -K_N^{n+\frac 12},v_N) -(R_N^{n+\frac 12},v_N) -(T_N^{n+\frac 12},v_N)+	\eps(v_N).
		\end{split}
	\end{equation}
	We also derive from   \eqref{spectral:sche2} that
	\begin{equation}\label{error:02}
		\frac{\hat e_N^{n+1}(\bz)-\tilde{e}_N^{n+1}(\bz)}{\delta t}=s_N^{n+1},\quad \forall \bz \in \Sigma_N,
			\end{equation}
		where   $s_N^{n+1}=-\lambda_N^{n+1}g'(u_N^{n+1})$.
	Denoting   $Q_N^{n+\frac 12}= \mathcal{N}(u(t^{n+\frac 12}))-\mathcal{N}(\frac 32\Pi_N u(t^{n})-\frac 12\Pi_N u(t^{n-1}))$, 
	we have 
	\begin{equation}
		\begin{split}
			&(\mathcal{N}(u(t^{n+\frac 12}))-\mathcal{N}(\frac 32u_N^n-\frac 12 u_N^{n-1}),v_N) =(Q_N^{n+\frac 12},v_N) \\&+(\mathcal{N}(\frac 32\Pi_N u(t^{n})-\frac 12\Pi_N u(t^{n-1}))-\mathcal{N}(\frac 32u_N^n-\frac 12 u_N^{n-1}),v_N) .
		\end{split}
	\end{equation}
		Then \eqref{error:01}  can be written as
	\begin{equation}\label{error:1}
		\begin{split}
			&(\frac{\tilde{e}_N^{n+1}-\hat e_N^n}{\delta t},v_N ) _N
			+(\nabla \frac{3\tilde e_N^{n+1}+\tilde e_N^{n-1}}{4},\nabla v_N)\\&+(\mathcal{N}(\frac 32\Pi_N u(t^{n})-\frac 12\Pi_N u(t^{n-1}))-\mathcal{N}(\frac 32u_N^n-\frac 12 u_N^{n-1}),v_N) \\&= -(K_N^{n+\frac 12},v_N) -(R_N^{n+\frac 12},v_N) -(T_N^{n+\frac 12},v_N) -(Q_N^{n+\frac 12},v_N)+	\eps(v_N).
		\end{split}
	\end{equation}
	Taking $v_N=2\delta t\tilde{e}_N^{n+1}$ in \eqref{error:1}, we obtain 
	\begin{equation}\label{error:3}
		\begin{split}
			&(\tilde{e}_N^{n+1}-\hat e_N^n,2\tilde{e}_N^{n+1})_N  + (R_N^{n+\frac 12}+K_N^{n+\frac 12}+T_N^{n+\frac 12}, 2\delta t\tilde{e}_N^{n+1}) +2\delta t(\nabla\frac{3\tilde{e}_N^{n+1}+\tilde e_N^{n-1}}{4},\nabla \tilde{e}_N^{n+1})
			\\&+(\mathcal{N}(\frac 32\Pi_N u(t^{n})-\frac 12\Pi_N u(t^{n-1}))-\mathcal{N}(\frac 32u_N^n-\frac 12 u_N^{n-1}), 2\delta t\tilde e_N^{n+1}) +(Q_N^{n+\frac 12},2\delta t\tilde e_N^{n+1})=2\delta t	\eps(\tilde e_N^{n+1}).
		\end{split}
	\end{equation}
		For the first term in \eqref{error:3}, we have
	\begin{equation}\label{error:4}
		\begin{split}
			&(\tilde{e}_N^{n+1}-\hat e_N^n,2\tilde{e}_N^{n+1})_N =\|\tilde{e}_N^{n+1}\|_N^2-\|\hat e_N^n\|_N^2 +\|\tilde{e}_N^{n+1}-\hat e_N^n\|_N^2.
		\end{split}
	\end{equation}
	We rewrite \eqref{error:02} as
	\begin{equation}\label{erro:eq:2}
		\begin{split}
			\hat e_N^{n+1}(\bz)-\delta t s_N^{n+1}(\bz)=\tilde{e}_N^{n+1}(\bz),\quad \forall \bz \in \Sigma_N,
		\end{split}
	\end{equation}
	and take the discrete inner product of \eqref{erro:eq:2} with itself to get
	\begin{equation}\label{generic:error:1}
		\begin{split}
			\|\hat e_N^{n+1}\|_N^2 +\delta t^2\|s_N^{n+1}\|_N^2-2\delta t(\hat e_N^{n+1},s_N^{n+1})_N =\|\tilde{e}_N^{n+1}\|_N^2.
		\end{split}
	\end{equation}
On the other hand,
	\begin{equation*}
		\begin{split}
		&	2\delta t(\nabla\frac{3\tilde{e}_N^{n+1}+\tilde e_N^{n-1}}{4}, \nabla\tilde{e}_N^{n+1}) =\frac{\delta t}{4}\{5(\nabla\tilde{e}_N^{n+1},\nabla\tilde e_N^{n+1}) -(\nabla\tilde{e}_N^{n-1},\nabla\tilde e_N^{n-1}) \\&+(\nabla(\tilde{e}_N^{n+1}+\tilde e^{n-1}_N),\nabla (\tilde e_N^{n+1}+\tilde e_N^{n-1})) \}.
		\end{split}
	\end{equation*}
		Combining the  above equations, we obtain
	\begin{equation}\label{error:k1}
		\begin{split}
		&	\|\hat e_N^{n+1}\|_N^2-\|\hat e_N^n\|_N^2+\|\tilde{e}_N^{n+1}-\hat e_N^n\|_N^2+\delta t^2\|s_{N}^{n+1}\|_N^2-2\delta t(\hat e_N^{n+1},s_{h}^{n+1})_N \\&+\frac{\delta t}{4}\{5(\nabla \tilde{e}_N^{n+1},\tilde e_N^{n+1}) -(\nabla \tilde{e}_N^{n-1},\nabla \tilde e_N^{n-1}) +(\nabla (\tilde{e}_N^{n+1}+\tilde e^{n-1}_N), \nabla(\tilde e_N^{n+1}+\tilde e_N^{n-1})) \}\\
			&=-(R_N^{n+\frac 12}+K_N^{n+\frac 12}+T_N^{n+\frac 12}+Q_N^{n+\frac 12}, 2\delta t\tilde{e}_N^{n+1}) -(\mathcal{N}(\frac 32\Pi_N u(t^{n})-\frac 12\Pi_N u(t^{n-1}))\\&-\mathcal{N}(\frac 32u_N^n-\frac 12 u_N^{n-1}), 2\delta t\tilde e_N^{n+1}) +2\delta t\eps(\tilde e_N^{n+1}).
		\end{split}
	\end{equation}
We now bound the terms on the righthand side as follows.

Firstly, consider the final term in \eqref{error:k1}, using \eqref{errdisinner},  we obtain
\begin{equation*}
\begin{split}
2\delta t\eps(\tilde e_N^{n+1}) &\leq 2C\delta t N^{-l}\|\tilde e_N^{n+1}\| \leq  2C\delta t N^{-l}\|\tilde e_N^{n+1}-\hat e_N^n\| +  2C\delta t N^{-l}\|\hat e_N^{n}\|\\&\leq 8C^2\delta t^2 N^{-2l} +\frac 18\|\tilde e_N^{n+1}-\hat e_N^n\|^2 + \delta t\|\hat e_N^n\|^2 + C^2\delta t N^{-2l}.
\end{split}
\end{equation*}

Thanks to  the KKT-condition $\lambda_N^{n+1}\ge 0$ and  $a\leq \Pi_N u(t^{n+1}) \leq b$, we find
	\begin{equation*}
		\begin{split}
			-2\delta t(\hat e_N^{n+1},s_{h}^{n+1})_N &= -2\delta t(u_N^{n+1}- \Pi_N u(t^{n+1}),\lambda_N^{n+1}g'(u_N^{n+1})) _N+2\delta t(\lambda_N^{n+1}, g(u_N^{n+1})) _N
			\\&= -2\delta t (\lambda_N^{n+1},-(u_N^{n+1})^2+ab+2\Pi_N u(t^{n+1})u_N^{n+1}-(a+b)\Pi_N u(t^{n+1}) )_N 
			\\&=2\delta t (\lambda_N^{n+1},(\Pi_N u(t^{n+1})-u_N^{n+1})^2)_N \\&-2\delta t ( \lambda_N^{n+1},(\Pi_N u(t^{n+1})-a)(\Pi_N u(t^{n+1})-b))_N \ge 0.
		\end{split}
	\end{equation*}
	On the other hand, since $ \mathcal{N}(u)=f_1(u)+\nabla\cdot f_2(u)$ with \eqref{assump1} and  \eqref{assump0},  we have
	\begin{equation}\label{non:error}
		\begin{split}
			&(\mathcal{N}(\frac 32\Pi_N u(t^{n})-\frac 12\Pi_N u(t^{n-1}))-\mathcal{N}(\frac 32u_N^n-\frac 12 u_N^{n-1}), 2\delta t\tilde e_N^{n+1}) 
			\\&=(f_1(\frac 32\Pi_N u(t^{n})-\frac 12\Pi_N u(t^{n-1}))-f_1(\frac 32u_N^n-\frac 12 u_N^{n-1}), 2\delta t\tilde e_N^{n+1})  \\&-(f_2(\frac 32\Pi_N u(t^{n})-\frac 12\Pi_N u(t^{n-1}))-f_2(\frac 32u_N^n-\frac 12 u_N^{n-1}), 2\delta t\Grad\tilde e_N^{n+1}) .
		\end{split}
	\end{equation}
The  terms on the righthand side of \eqref{non:error} can be bounded as follows:
	\begin{equation}\label{non:error:2}
		\begin{split}
			&(f_1(\frac 32\Pi_N u(t^{n})-\frac 12\Pi_N u(t^{n-1}))-f_1(\frac 32u_N^n-\frac 12 u_N^{n-1}), 2\delta t\tilde e_N^{n+1})  
			\leq 2C_1\delta t(|\frac 32 \hat e_N^n-\frac 12 \hat  e_N^{n-1}|, \tilde e_N^{n+1}) \\& =2C_1\delta t(|\frac 32 \hat e_N^n-\frac 12 \hat  e_N^{n-1}|, \tilde e_N^{n+1}-\hat e_N^n)  + 2C_1\delta t(|\frac 32 \hat e_N^n-\frac 12 \hat  e_N^{n-1}|,\hat e_N^n) \\
	&  \leq \frac 18\|\tilde e_N^{n+1}-\hat e_N^n\|^2 +8C_1^2\delta t^2\|\frac 32 \hat e_N^n-\frac 12 \hat  e_N^{n-1}\|^2
		+ C_1 \delta t(\|\hat e_N^n\|^2+\|\frac 32 \hat e_N^n-\frac 12 \hat  e_N^{n-1}\|^2).
		\end{split}
	\end{equation}
		Similarly, 
	\begin{equation}
		\begin{split}
			&(f_2(\frac 32\Pi_N u(t^{n})-\frac 12\Pi_N u(t^{n-1}))-f_2(\frac 32u_N^n-\frac 12 u_N^{n-1}), 2\delta t\Grad\tilde e_N^{n+1}) \\&
			\leq 2\delta t C_2(|\frac 32 \hat e_N^n-\frac 12 \hat  e_N^{n-1}|, \Grad \tilde e_N^{n+1}) \leq  \frac 13\delta t \|\Grad \tilde e_N^{n+1}\|^2+3\delta tC_2^2\|\frac 32 \hat e_N^n-\frac 12 \hat  e_N^{n-1}\|^2.
		\end{split}
	\end{equation}
	It remains to deal with the first term on the righthand side of \eqref{error:k1}.
	\begin{equation}
		\begin{split}
			&-2\delta t(R_N^{n+\frac 12}+T_N^{n+\frac 12},\tilde{e}_N^{n+1}) =-2\delta t(R_N^{n+\frac 12}+T_N^{n+\frac 12},\tilde{e}_N^{n+1}-\hat e_N^n)  - 2\delta t(R_N^{n+\frac 12}+T_N^{n+\frac 12},\hat e_N^n) 
			\\&\leq 4\delta t^2\|R_N^{n+\frac 12}\|^2+4\delta t^2\|T_N^{n+\frac 12}\|^2+\frac 14\|\tilde{e}_N^{n+1}-\hat e_N^n\|^2 +\delta t(\|R_N^{n+\frac 12}\|^2+\|T_N^{n+\frac 12}\|^2+\|\hat e_N^n\|^2);
		\end{split}
	\end{equation}
and
	\begin{equation}\label{project:er}
		\begin{split}
		&	-2\delta t(K_N^{n+1},\tilde e_N^{n+1}) =-2\delta t(\frac{\bar e_N^{n+1}-\bar e_N^n}{\delta t}, \tilde e_N^{n+1}) \\
		&	=-2((I-\Pi_N)(u(t^{n+1})-u(t^n)), \tilde e_N^{n+1}-\hat e_N^n+\hat e_N^n) 
			\\&\leq 2|((I-\Pi_N)(u(t^{n+1})-u(t^n)), \tilde e_N^{n+1}-\hat e_N^n ) |
			+2|((I-\Pi_N)(u(t^{n+1})-u(t^n)), \hat e_N^n )| \\
	& \le  8\delta t \int_{t^n}^{t^{n+1}}\|(I-\Pi_N)u_t(t)\|^2 dt +\frac 18\|\tilde e_N^{n+1}-\hat e_N^n\|^2
+ \int_{t^n}^{t^{n+1}}\|(I-\Pi_N)u_t(t)\|^2 dt +\delta t\|\hat e_N^n\|^2;
	\end{split}
	\end{equation}
and
	\begin{equation}\label{non:ac}
		\begin{split}
			&(Q_N^{n+\frac 12},2\delta t\tilde e_N^{n+1}) 
			=(\mathcal{N}(u(t^{n+\frac 12}))-\mathcal{N}(\frac 32\Pi_N u(t^{n})-\frac 12\Pi_N u(t^{n-1})), 2\delta t\tilde e_N^{n+1}) 
			\\&=(\mathcal{N}(u(t^{n+\frac 12}))-\mathcal{N}(\frac 32 u(t^{n})-\frac 12 u(t^{n-1})), 2\delta t\tilde e_N^{n+1}) \\&+ (\mathcal{N}(\frac 32 u(t^{n})-\frac 12 u(t^{n-1}))-\mathcal{N}(\frac 32\Pi_N u(t^{n})-\frac 12\Pi_N u(t^{n-1})), 2\delta t\tilde e_N^{n+1}) .
		\end{split}
	\end{equation}
	For the first term in right hand side of \eqref{non:ac}, we have
	\begin{equation}\label{non:split}
		\begin{split}
			&(\mathcal{N}(u(t^{n+\frac 12}))-\mathcal{N}(\frac 32 u(t^{n})-\frac 12 u(t^{n-1})), 2\delta t\tilde e_N^{n+1}) 
			\\&=(f_1(u(t^{n+\frac 12}))-f_1(\frac 32u(t^n)-\frac 12 u(t^{n-1})), 2\delta t\tilde e_N^{n+1})  \\&-(f_2(u(t^{n+\frac 12}))-f_2(\frac 32u(t^n)-\frac 12 u(t^{n-1})), 2\delta t\Grad\tilde e_N^{n+1}) .
		\end{split}
	\end{equation}
	Using assumptions  \eqref{assump1}-\eqref{assump0} and Young's inequality, we have 
	\begin{equation}\label{non:split:2}
		\begin{split}
			&(f_1(u(t^{n+\frac 12}))-f_1(\frac 32u(t^n)-\frac 12 u(t^{n-1})), 2\delta t\tilde e_N^{n+1}) \leq C_1 (|J_N^{n+\frac 12}|,2\delta t\tilde e_N^{n+1})  \\& =
			2C_1\delta t (|J_N^{n+\frac 12}|,\tilde e_N^{n+1}-\hat e_N^n)  + 2C_1\delta t (|J_N^{n+\frac 12}|,\hat e_N^n) \\
	&  \leq \frac 14\|\tilde e_N^{n+1}-\hat e_N^n\|^2 +4C_1^2 \delta t^2 \|J_N^{n+\frac 12}\|^2+ \delta t C_1(\|\hat e_N^n\|^2+\|J_N^{n+\frac 12}\|^2);
		\end{split}
	\end{equation}
	and
	\begin{equation}\label{non:split:3}
		\begin{split}
			&(f_2(u(t^{n+\frac 12}))-f_2(\frac 32u(t^n)-\frac 12 u(t^{n-1})), 2\delta t\Grad\tilde e_N^{n+1})  \leq C_2(|J_N^{n+\frac 12}|,  2\delta t\Grad\tilde e_N^{n+1}) \\& \leq 3\delta tC_2^2\|J_N^{n+\frac 12}\|^2+\frac 13\delta t\|\Grad \tilde e_N^{n+1}\|^2.
		\end{split}
	\end{equation}	
For the second term in right hand side of  \eqref{non:ac}, similar with \eqref{non:split:2} and \eqref{non:split:3}, we have
\begin{equation}\label{non:split:4}
\begin{split}
&(\mathcal{N}(\frac 32 u(t^{n})-\frac 12 u(t^{n-1}))-\mathcal{N}(\frac 32\Pi_N u(t^{n})-\frac 12\Pi_N u(t^{n-1})), 2\delta t\tilde e_N^{n+1}) 
\\&\leq C_1(|\frac 32 \bar e_N^n-\frac 12 \bar e_N^{n-1}|, 2\delta t\tilde e_N^{n+1}) + C_2(|\frac 32 \bar e_N^n-\frac 12 \bar e_N^{n-1}|, 2\delta t\nabla \tilde e_N^{n+1}).
\end{split}
\end{equation}
For the first term in the right hand side of \eqref{non:split:4}, using assumption \eqref{assump3},  we have
\begin{equation*}
\begin{split}
&C_1(|\frac 32 \bar e_N^n-\frac 12 \bar e_N^{n-1}|, 2\delta t\tilde e_N^{n+1}) \leq 2C_1\delta t(|\frac 32 \bar e_N^n-\frac 12 \bar e_N^{n-1}|, \tilde e_N^{n+1}-\hat e_N^n) +2C_1\delta t(|\frac 32 \bar e_N^n-\frac 12 \bar e_N^{n-1}|, \hat e_N^{n})\\& \leq C \delta t N^{-2l} +\frac 18\|\tilde e_N^{n+1}-\hat e_N^n\|^2 + C_1\delta t \|\hat e_N^n\|^2.
\end{split}
\end{equation*}
For the second term in the right hand side of \eqref{non:split:4}, we have
\begin{equation*}
\begin{split}
& C_2(|\frac 32 \bar e_N^n-\frac 12 \bar e_N^{n-1}|, 2\delta t\nabla \tilde e_N^{n+1})\leq \frac 13\delta t\|\Grad \tilde e_N^{n+1}\|^2 +3\delta tC_2^2\|\frac 32 \bar e_N^n-\frac 12 \bar e_N^{n-1}\|^2 \\&\leq \frac 13\delta t\|\Grad \tilde e_N^{n+1}\|^2 + C\delta t N^{-2l}.
 \end{split}
\end{equation*}
Combining the above relations into  \eqref{error:k1} and using \eqref{disinner},  we arrive at 
	\begin{equation}\label{error:final:1}
		\begin{split}
			&\|\hat e_N^{n+1}\|_N^2-\|\hat e_N^n\|_N^2+\delta t^2\|s^{n+1}_{N}\|_N^2+\frac{\delta t}{4}\{(\nabla \tilde{e}_N^{n+1},\nabla\tilde e_N^{n+1}) -(\nabla \tilde{e}_N^{n-1},\nabla\tilde e_N^{n-1}) \\&+(\nabla (\tilde{e}_N^{n+1}+\tilde e^{n-1}_N),\nabla(\tilde e_N^{n+1}+\tilde e_N^{n-1})) \} \leq (8C_1^2\delta t^2+C_1\delta t+3\delta tC_2^2)\|\frac 32 \hat e_N^n-\frac 12 \hat  e_N^{n-1}\|^2\\&+3(C_1+1) \delta t\|\hat e_N^n\|_N^2
			+(4\delta t^2+\delta t)\|R_N^{n+\frac 12}\|^2+(4\delta t^2+\delta t)\|T_N^{n+\frac 12}\|^2 \\&
			+(4C_1^2 \delta t^2+\delta tC_1+3\delta t C_2^2) \|J_N^{n+\frac 12}\|^2+ 8\delta t \int_{t^n}^{t^{n+1}}\|(I-\Pi_N)u_t(t)\|^2 dt\\&+ \int_{t^n}^{t^{n+1}}\|(I-\Pi_N)u_t(t)\|^2 dt+(C^2+2C+8C^2\delta t)\delta t N^{-2l} ,\quad\forall n\ge 1 .
		\end{split}
	\end{equation}
For $n=0$, a similar estimate can be easily derived. 
	Summing up  \eqref{error:final:1} from $n=1$ to $n=m-1$ and its corresponding inequality at $n=0$, we obtain
	\begin{equation}\label{err:gr}
		\begin{split}
			&\|\hat e_N^{m}\|_N^2+\delta t^2\sum\limits_{n=1}^{m-1}\|s^{n+1}_{N}\|_N^2+\frac{\delta t}{4}\|\nabla \tilde{e}_N^{m}\|^2 +\frac{\delta t}{4}\|\nabla \tilde{e}_N^{m-1}\|^2+\sum\limits_{n=1}^{m-1}\|\nabla (\tilde{e}_N^{n+1}+\tilde e^{n-1}_N)\|^2 \\& \leq \|\hat e_N^0\|_N^2+\frac{\delta t}{4}\|\nabla \tilde{e}_N^{0}\|^2++\frac{\delta t}{4}\|\nabla \tilde{e}_N^{1}\|^2+\sum\limits_{n=0}^{m-1}\{ (8C_1^2\delta t^2+C_1\delta t+3\delta tC_2^2)\|\frac 32 \hat e_N^n-\frac 12 \hat  e_N^{n-1}\|^2\\&+3(C_1+1) \delta t\|\hat e_N^n\|^2
			+(4\delta t^2+\delta t)\|R_N^{n+\frac 12}\|^2+(4\delta t^2+\delta t)\|T_N^{n+\frac 12}\|^2 \\&
			+(4C_1^2 \delta t^2+\delta tC_1+3\delta t C_2^2) \|J_N^{n+\frac 12}\|^2\}+ 8\delta t \int_{0}^{T}\|(I-\Pi_N)u_t(t)\|^2 dt\\&+ \int_{0}^{T}\|(I-\Pi_N)u_t(t)\|^2 dt+(C^2+2C+8C^2\delta t)T  N^{-2l}.
		\end{split}
	\end{equation}
	For the term in \eqref{err:gr} with $l\ge 0$, we have
	\begin{equation*}
		\begin{split}
			&\int_{0}^{T}\|(I-\Pi_N)u_t(t)\|^2 dt  \leq CN^{-2l} \|u_t\|_{L^2(0,T;H^l)}^2;\quad 
			\|T_N^{n+\frac 12}\|^2  \leq  C\delta t^3 \int_{t^n}^{t^{n+1}}\|u_{tt}\|^2_{H^2}dt;\\
			&\|J_N^{n+\frac 12}\|^2  \leq  C\delta t^3\int_{t^n}^{t^{n+1}}\|u_{tt}\|^2dt;\quad 
			\|R_N^{n+\frac 12}\|^2  \leq C\delta t^3\int_{t^n}^{t^{n+1}}\|u_{ttt}\|^2dt.
		\end{split}
	\end{equation*}
	Finally, applying the discrete Gronwall's Lemma to the above, using the norm equivalence \eqref{disinner} and the triangular inequality, we obtain the desired  result.
\end{proof}

\begin{remark}
	By following exactly the same procedure, we can also derive a similar error estimate if we use a hybrid Fourier spectral method instead of the hybrid Legendre spectral method.
\end{remark}
\section{Some typical applications}
The bound preserving schemes that we constructed and studied in previous sections can be applied to a large class of PDEs which are bound preserving. We describe applications to several typical examples below.
\subsection{Allen-Cahn equation}
Consider the Allen-Cahn equation \cite{allen1979microscopic}  
\begin{equation}\label{allen:cahn}
u_t- \Delta u +\frac{1}{\eps^2}u(u^2-1)=0,
\end{equation}
with homogeneous Dirichlet, homogeneous Neumann or periodic boundary condition, and $\eps$ is a positive constant. It is well known that the above equation satisfies the maximum principle, in particular, if the values of the initial condition $u_0$ is in $[-1,1]$,  the solution of the Allen-Cahn equation \eqref{allen:cahn} will stay within the range $[-1,1]$.  
Setting $\mathcal{L}=-\Delta +\frac{1}{\eps^2}$ and $\mathcal{N}(u)=f_1(u)=\frac{1}{\eps^2}u(u^2-1)-\frac{1}{\eps^2} $, a second-order scheme based on the modified Crank-Nicholson for \eqref{allen:cahn} is:
\begin{equation}\label{semi:allen}
	\begin{split}
		&\frac{\tilde u^{n+1}-u^n}{\delta t} +\mathcal{L}(\frac 34 \tilde u^{n+1}+\frac 14 \tilde u^{n-1})+\mathcal{N}(\frac 32u^n-\frac 12u^{n-1})=\lambda^ng'(u^n);
	\end{split}
\end{equation}
and
\begin{equation}\label{semi:allen2}
	\begin{split}		
		&\frac{u^{n+1}-\tilde u^{n+1}}{\delta t}=\frac12(\lambda^{n+1}g'(u^{n+1})-\lambda^ng'(u^n)),\\
		&\lambda^{n+1}\ge 0, \; g(u^{n+1}) \ge 0,\;\lambda^{n+1}g(u^{n+1})=0,
	\end{split}
\end{equation}
where $g(u)=(1+u)(1-u)$.

Similarly, we have its cut off version:
\begin{equation}\label{cutoff:allen}
	\begin{split}
		&\frac{\tilde u^{n+1}-u^n}{\delta t} +\mathcal{L}(\frac 34 \tilde u^{n+1}+\frac 14 \tilde u^{n-1})+\mathcal{N}(\frac 32u^n-\frac 12u^{n-1})=0;
	\end{split}
\end{equation}
and
\begin{equation}\label{cutoff:allen2}
	\begin{split}		
		&\frac{u^{n+1}-\tilde u^{n+1}}{\delta t}=\lambda^{n+1}g'(u^{n+1}),\\
		&\lambda^{n+1}\ge 0, \; g(u^{n+1}) \ge 0,\;\lambda^{n+1}g(u^{n+1})=0.
	\end{split}
\end{equation}

Since $f_1(u)=0$, and $f_2(u)$ is certainly locally Lipschitz and satisfies \eqref{assump1}- \eqref{assump2}.
 Hence, results which are similar to those in Theorem 3.1 and  Theorem 4.1 can be derived  for the above schemes.

\subsection{Cahn-Hilliard equation with variable mobility}
Consider the  Cahn-Hilliard equation \cite{cahn1958free}  with a logarithmic  potential: 
\begin{equation}\label{CH}
\begin{split}
& u_t =\Grad\cdot(M(u)\Grad\mu) ,\\
&\mu= -\eps^2\Delta u +\ln(1+u)-\ln(1-u)-\theta_0u,
\end{split}
\end{equation}
where $\mu$ is the chemical potential and $M(u)=1-u^2>0$ is the mobility function. $\theta_0,\eps$  are two positive constants. $u$ and $\mu$ are prescribed with homogeneous Neumann or periodic  boundary condition. The Cahn-Hilliard equation \eqref{CH}  is a gradient flow which takes on the form
\begin{equation}
	u_t=\nabla \cdot(M(u)\Grad \frac{\delta E}{\delta u}),
\end{equation}
 with  the total free energy 
\begin{equation}
E(u)=\int_{\Omega} (1+u)\ln(1+u) +(1-u)\ln(1-u)-\frac{\theta_0}{2}u^2+\frac{\eps^2}{2}|\Grad u|^2 d\bx.
\end{equation}
With a given initial condition $\|u_0\|_{L^\infty} < 1-\gamma$ for a constant $\gamma\in (0,1)$, due to the singular logarithmic  potential, the  solution of Cahn-Hilliard equation \eqref{CH} is expected to remain in the range  $(-1+\delta,1-\delta)$ for some $\delta\in (0,1)$  \cite{debussche1995cahn,elliott1996cahn}.  Note that   \eqref{CH} is a fourth-order equation written as a system of two coupled second-order equations, so the approach for constructing bound preserving schemes introduced in Section 2 can be directly applied to   \eqref{CH}. For example, the second-order version of \eqref{high:bound:lag:1}-\eqref{high:bound:lag:2} for   \eqref{CH} is as follows:

\begin{equation}\label{bdf2:CH:LM}
\begin{split}
& \frac{3\tilde{u}^{n+1}-4u^n+u^{n-1}}{2\delta t} =\Grad\cdot(M(2u^n-u^{n-1})\Grad\mu^{n+1}) +\lambda^{n} g'(u^n),\\
&\mu^{n+1}= -\eps^2\Delta u^{n+1} +\ln(1+2u^n-u^{n-1})-\ln(1-2u^n+u^{n-1})-\theta_0(2u^n-u^{n-1});
\end{split}
\end{equation}
and
\begin{equation}\label{bdf2:CH:LM2}
\begin{split}
&\frac{3u^{n+1}-3\tilde{u}^{n+1}}{2\delta t}=\lambda^{n+1}g'(u^{n+1})-\lambda^{n} g'(u^n),\\
&\lambda^{n+1} \ge 0,\; g(u^{n+1}) \ge 0,\;\lambda^{n+1} g(u^{n+1})=0,
\end{split}
\end{equation}
where $g(u)=(u+1-\delta)(1-\delta-u)$.  Notice that $g(u)=(u+1-\delta)(1-\delta-u)>0$ is equivalent to $-1+\delta \leq u \leq 1-\delta$.
 
The system \eqref{CH} also preserves mass. Indeed, integrate the first equation in \eqref{CH} over $\Omega$, we obtain  $\partial_t \int_\Omega u d\bx =0$.  As described in Section 2, we can also easily modify the scheme \eqref{bdf2:CH:LM}-\eqref{bdf2:CH:LM2} to construct a bound and mass preserving scheme for \eqref{CH}.

While the stability results in Section 3 was derived only for a second-order equation for the sake of simplicity, since the nonlinear term   $\mathcal{N}(u)=f_2(u)=\ln(1+u)-\ln(1-u)-\theta_0u$ is locally Lipschitz for $u\in (-1,1)$ and satisfies \eqref{assump1}-\eqref{assump2}, a similar procedure can be used to derive a stability result which is similar to Theorem 3.1. However, the error analysis  in Section 4 can not be easily extended to this case.

\subsection{Fokker-Planck equation}
Consider the following Fokker-Planck equation  
\begin{equation}\label{fockker-planck}
\partial_t u = \partial_x(x u(1-u)+\partial_x u),
\end{equation}
with no flux or periodic boundary conditions, which models the relaxation of fermion and boson gases  taking on the form \cite{carrillo20081d,sun2018discontinuous}. The long time asymptotics of the one dimensional model has been studied in \cite{carrillo20081d}.  

The Fokker-Planck equation \eqref{fockker-planck}  can be interpreted as a gradient flow
\begin{equation}
	\partial_t u = \partial_x (u(1-u)\partial_x \frac{\delta E}{\delta u}),
\end{equation}
with $E(u)$ being the entropy functional
\begin{equation}
E(u) =\int_{\Omega} \big( \frac{x^2}{2}u + u\log(u)+(1-u)\log(1-u) \big) d\bx.
\end{equation}
Hence, the solution of \eqref{fockker-planck}  is expected to take values  in $[0,1]$.

The approach for constructing bound preserving schemes introduced in Section 2 can be directly applied to   \eqref{fockker-planck}. For example, let $\mathcal{L} u=-\partial_{xx} u$ and $\mathcal{N}(u)=\partial_x f_2(u)=\partial_x(-xu (1-u)) $, a second-order version of \eqref{high:bound:lag:1}-\eqref{high:bound:lag:2} for   \eqref{fockker-planck} is as follows: 
\begin{equation}\label{scheme:fockker:planck}
\begin{split}
&\frac{3u^{n+1}-4u^{n}+u^{n-1}}{2\delta t} = \partial_x(x(2u^n-u^{n-1})(1-2u^n+u^{n-1})+\partial_x u^{n+1}) + \lambda^n g'(u^n);
\end{split}
\end{equation}
and
\begin{equation}\label{scheme:fockker:planck2}
	\begin{split}
&\frac{3\tilde u^{n+1}-3\tilde u^{n+1}}{2\delta t} = \lambda^{n+1} g'(u^{n+1})-\lambda^n g'(u^n),\\
&\lambda^{n+1} \ge 0,\; g(u^{n+1}) \ge 0,\;\lambda^{n+1} g(u^{n+1})=0,
\end{split}
\end{equation}
where $g(u)=u(1-u)$. 

We observe that the Fokker-Plank equation  \eqref{fockker-planck} with  no flux or periodic boundary conditions conserves mass, i.e.,   $\partial_t \int_\Omega u d\bx =0$. The above scheme can  be easily modified to be  mass conserving as follows:
\begin{equation}\label{mass:fockker:planck}
	\begin{split}
		&\frac{3u^{n+1}-4u^{n}+u^{n-1}}{2\delta t} = \partial_x(x(2u^n-u^{n-1})(1-2u^n+u^{n-1})+\partial_x u^{n+1}) + \lambda^n g'(u^n),
	\end{split}
\end{equation}
and
\begin{equation}\label{mass:fockker:planck2}
	\begin{split}		
		&\frac{3u^{n+1}-3\tilde u^{n+1}}{2\delta t} = \lambda^{n+1} g'(u^{n+1})-\lambda^n g'(u^n)+\xi^{n+1},\\
		&\lambda^{n+1} \ge 0,\; g(u^{n+1}) \ge 0,\;\lambda^{n+1} g(u^{n+1})=0,\;(u^{n+1},1)=(u^n,1).
	\end{split}
\end{equation}

It is clear that $f_2(u)=-xu (1-u)$ is locally Lipschitz and satisfies \eqref{assump1}- \eqref{assump0} with $f_1(u)=0$. Therefore,  a similar result as  in Theorem 3.1 can be derived for the scheme \eqref{scheme:fockker:planck}-\eqref{scheme:fockker:planck2} and \eqref{mass:fockker:planck}-\eqref{mass:fockker:planck2}.

\section{Numerical results}
In this section, we will present various numerical experiments to validate the proposed bound preserving schemes. For all examples presented below, we assume periodic boundary conditions in $\Omega=[0,2\pi)^d$, and use 
a Fourier-spectral  method for spatial approximation.

\subsection{Allen-Cahn equation}
The first example is the Allen-Cahn equation \eqref{allen:cahn}.
\subsubsection{Accuracy test}
We first verify the convergence rate for the scheme  \eqref{semi:allen}-\eqref{semi:allen2}  and its first-order version for \eqref{full:dis} in the domain $\Omega=[0,2\pi]^2$
with the initial condition 
\begin{equation}\label{ini:allen}
u(x,y,0)=\tanh(\frac{1-\sqrt{(x-\pi)^2+(y-\pi)^2}}{\sqrt{2}\eps}).
\end{equation}
 We use  $128^2$ uniform collocation points in $[0,2\pi]^2$, i.e.,  $\Sigma_N=\{x_{jk}=(\frac{j}{2\pi},\frac{k}{2\pi}); j,k=0,1,,\cdots,128\}$, so that the spatial discretization error is negligible compared with the time discretization error. We shall test their accuracy as approximations of \eqref{full:dis} and  \eqref{strong} respectively.
 
\begin{figure}[htbp]
\centering
\includegraphics[width=0.45\textwidth,clip==]{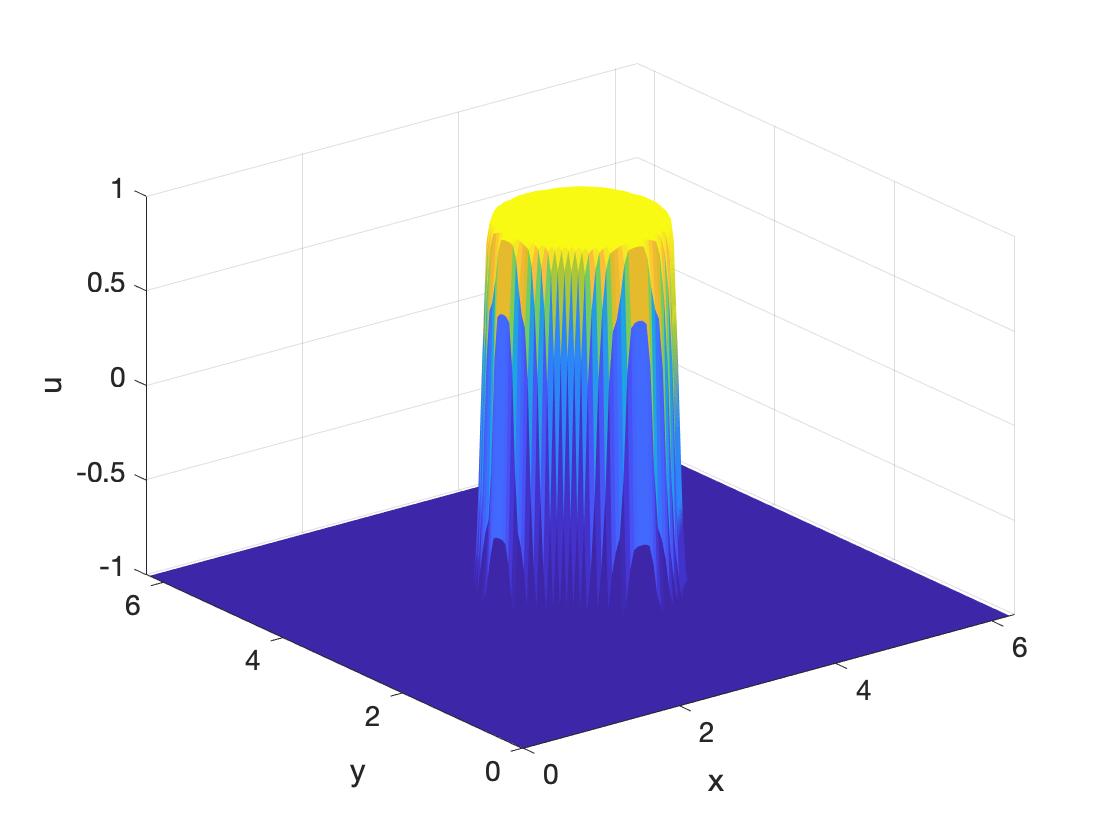}
\includegraphics[width=0.45\textwidth,clip==]{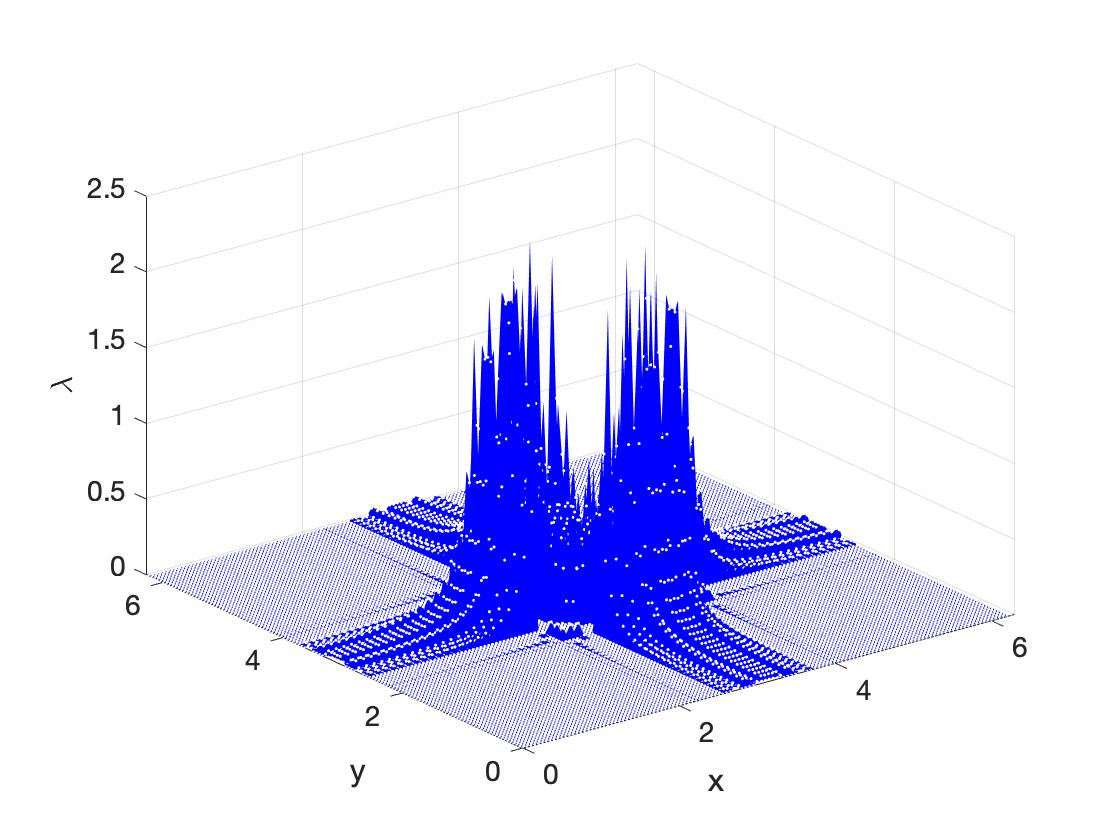}
\caption{ Numerical solution  $u$ and  Lagrange multiplier $\lambda$ at $T=0.001$ computed by  scheme \eqref{high:bound:lag:1}-\eqref{high:bound:lag:2} with  $k=2$ and $\delta t= 10^{-6}$.}\label{ex_1}
\end{figure}

  First, we consider  these schemes as  approximations of \eqref{full:dis}, and use  the reference solution   computed by  \eqref{semi:allen}-\eqref{semi:allen2}  with a very small time step $\delta t=10^{-6}$.
 We observe from table \ref{table1} that the scheme  \eqref{semi:allen}-\eqref{semi:allen2} (resp. its first-order version)  achieves  second-order (resp. first-order)  convergence rate in time.  The scheme \eqref{cutoff:allen}-\eqref{cutoff:allen2} only achieves the first-order convergence in time.  We plot in  Fig.\;\ref{ex_1} the profile of numerical solution $u$  and the   Lagrange multiplier $\lambda$ at $T=0.001$.

  \begin{table}[ht]
  	\centering
  	\begin{tabular}{r||c|c|c|c|c|c|}
  		\hline
  		$\delta t$       &     BDF1 version of   \eqref{semi:allen}-\eqref{semi:allen2} & Order &   \eqref{semi:allen}-\eqref{semi:allen2}  & Order  & \eqref{cutoff:allen}-\eqref{cutoff:allen2} & Order \\ \hline
  		$4\times 10^{-5}$    &$4.89E(-3)$  & $-$  &$3.56E(-4)$ &$-$   & $1.36E(-3) $ & $-$    \\\hline
  		$2\times 10^{-5}$     &$2.47E(-3)$ & $0.98$&$9.50E(-5)$ &$1.90$ & $6.75E(-4)$ & $1.01$    \\\hline
  		$1\times 10^{-5}$     &$1.24E(-3)$ &$0.99$  &$2.31E(-5)$ &$2.04$ & $3.24E(-4)$ & $1.06$  \\\hline
  		$5\times 10^{-6}$  &$6.22E(-4)$ &$0.99$ &$5.84E(-6)$ &$1.98$  & $1.44E(-4)$ & $1.17$ \\ \hline
  		$2.5\times 10^{-6}$  &$3.11E(-4)$&$1.00$ &$1.25E(-6)$ &$2.22$ & $5.43E(-5)$ & $1.40$ \\\hline
  		\hline
  	\end{tabular}
  	\vskip 0.5cm
  	\caption{Accuracy test for approximations to \eqref{full:dis}: The $L^{\infty}$ errors   at $t=0.01$   with  $\eps^2=0.001$.}\label{table1}
  \end{table}

    \begin{table}[ht]
  	\centering
  	\begin{tabular}{r|c|c|c|c|c|c}
  		\hline
  		$\delta t$       &   \eqref{semi:allen}-\eqref{semi:allen2}  & Order  & \eqref{cutoff:allen}-\eqref{cutoff:allen2} & Order \\ \hline
  		$4\times 10^{-5}$   &$1.05E(-4)$ &$-$   & $1.05E(-4) $ & $-$    \\\hline
  		$2\times 10^{-5}$    &$4.25E(-5)$ &$1.30$ & $4.25E(-5)$ & $1.30$    \\\hline
  		$1\times 10^{-5}$    &$1.00E(-5)$ &$2.08$ & $1.00E(-5)$ & $2.08$  \\\hline
  		$5\times 10^{-6}$ &$2.76E(-6)$ &$1.86$  & $2.76E(-6)$ & $1.86$ \\ \hline
  		$2.5\times 10^{-6}$ &$6.29E(-7)$ &$2.13$ & $6.29E(-7)$ & $2.13$ \\\hline
  		\hline
  	\end{tabular}
  	\vskip 0.5cm
  	\caption{Accuracy test for approximations to \eqref{strong}: The $L^{\infty}$ errors   at $t=0.01$   with  $\eps^2=0.001$.}\label{table2}
  \end{table}
  
 We then consider  these schemes as  approximations of \eqref{strong}, and use the reference solution as a highly accurate approximation to the original PDE \eqref{strong}  which is computed by a standard semi-implicit scheme with $\delta t=10^{-8}$. We compare the accuracy between the scheme  \eqref{semi:allen}-\eqref{semi:allen2}  and its cur-off version  \eqref{cutoff:allen}-\eqref{cutoff:allen2}. The results are reported in Table \ref{table2}. We observe that both schemes have essentially the same accuracy and are  second-order in time, which are consistent with the error estimates in Theorem 4.1.
   
   The results reported in Tables 1 and 2 are consistent with {\bf Remark 2.1}.
 
\subsubsection{Comparison with a usual semi-implicit scheme}
We consider the  Allen-Cahn equation with  $\eps^2=0.001$ and the initial condition 
\begin{equation}\label{ini:allen:2}
\begin{split}
u(x,y,0)&=\tanh(\frac{1-\sqrt{(x-\pi)^2+(y-3\pi/2)^2}}{\sqrt{2}\eps})
\\&+\tanh(\frac{1-\sqrt{(x-\pi)^2+(y-3\pi/4)^2}}{\sqrt{2}\eps})+1.
\end{split}
\end{equation}
 We use the scheme \eqref{cutoff:allen}-\eqref{cutoff:allen2}  and its usual semi-implicit version:
 \begin{equation}\label{cutoff:allenB}
 	\begin{split}
 		&\frac{ u^{n+1}-u^n}{\delta t} +\mathcal{L}(\frac 34 u^{n+1}+\frac 14  u^{n-1})+\mathcal{N}(\frac 32u^n-\frac 12u^{n-1})=0,
 	\end{split}
 \end{equation} 
  with time step  $\delta t=8\times 10^{-4}$ and  $128^2$ Fourier modes.
 
\begin{figure}[htbp]
\centering
\subfigure[$u$ at $T=0.08 $ by \eqref{cutoff:allenB}.]{
\includegraphics[width=0.30\textwidth,clip==]{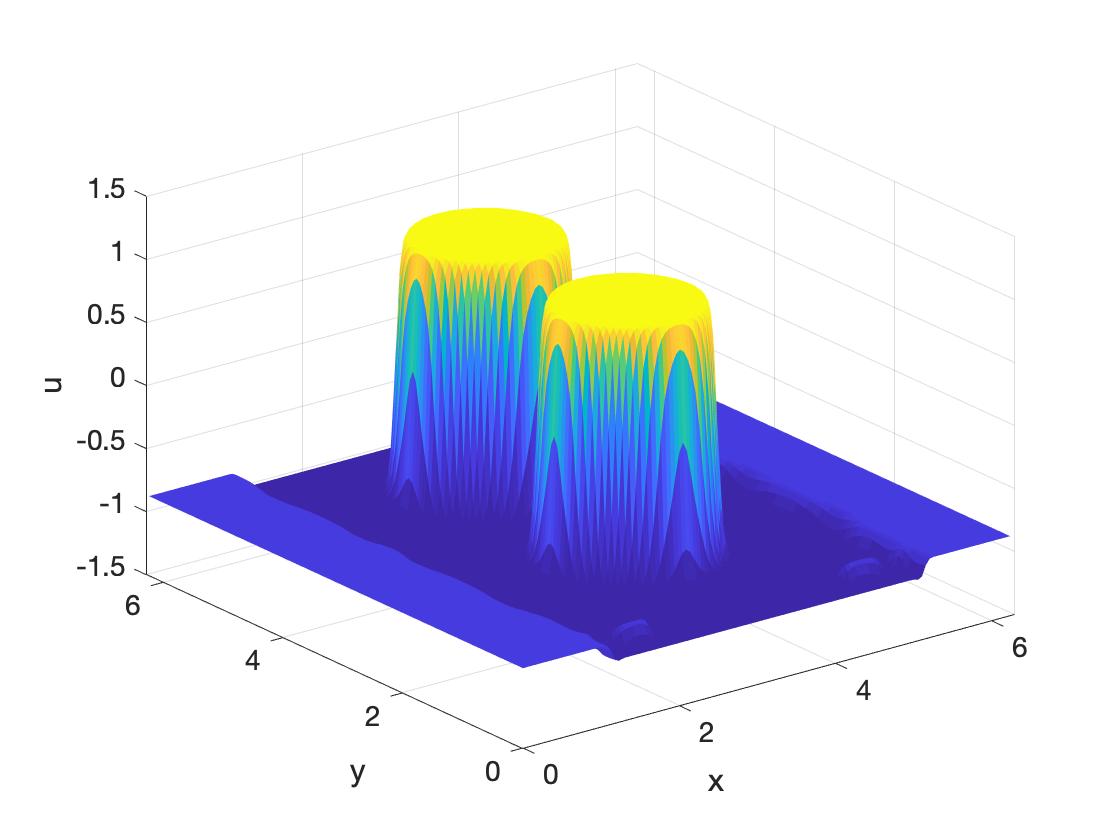}}
\subfigure[$u$ at $T=0.08$ by \eqref{cutoff:allen}-\eqref{cutoff:allen2}.]{
\includegraphics[width=0.30\textwidth,clip==]{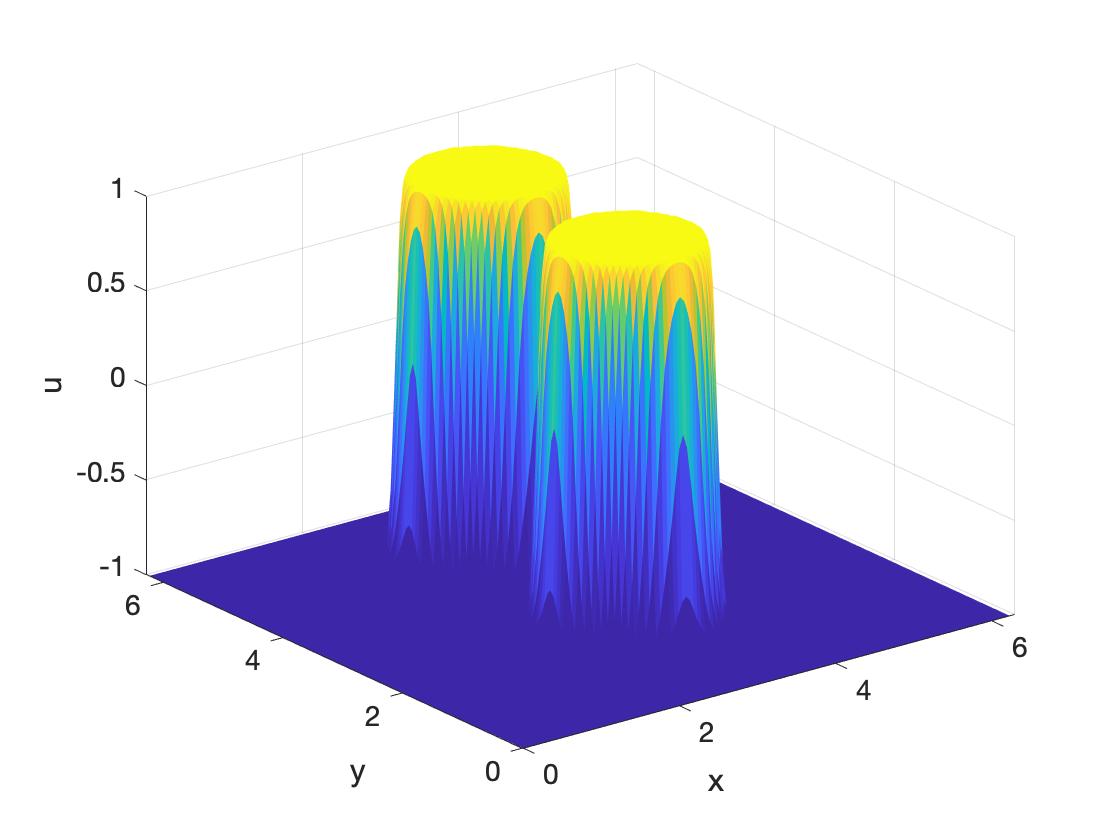}}
\subfigure[$\lambda$ at $T=0.08$.]{
\includegraphics[width=0.30\textwidth,clip==]{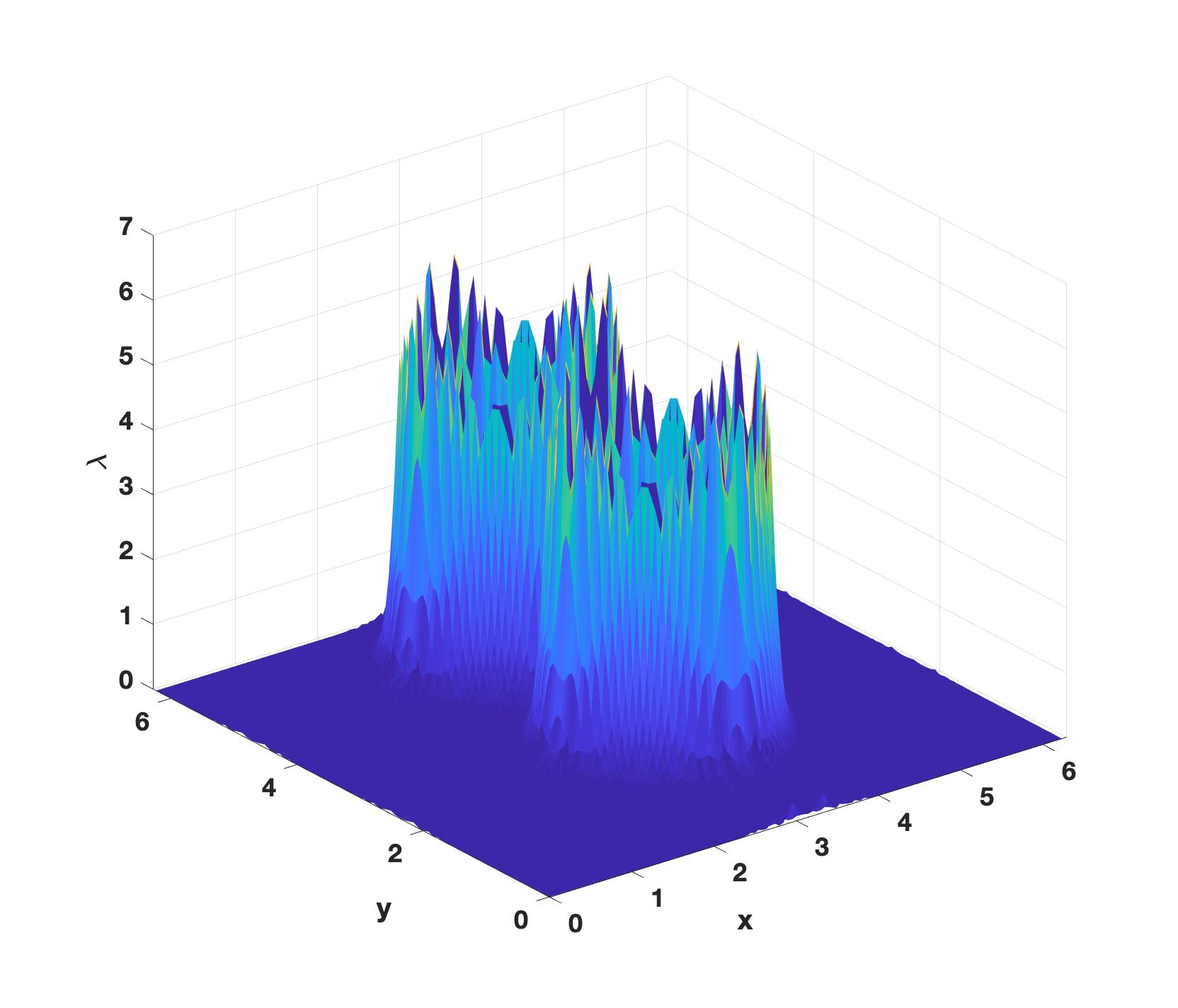}}
\subfigure[$u$ at $T=0.4$ by \eqref{cutoff:allenB}.]{
\includegraphics[width=0.30\textwidth,clip==]{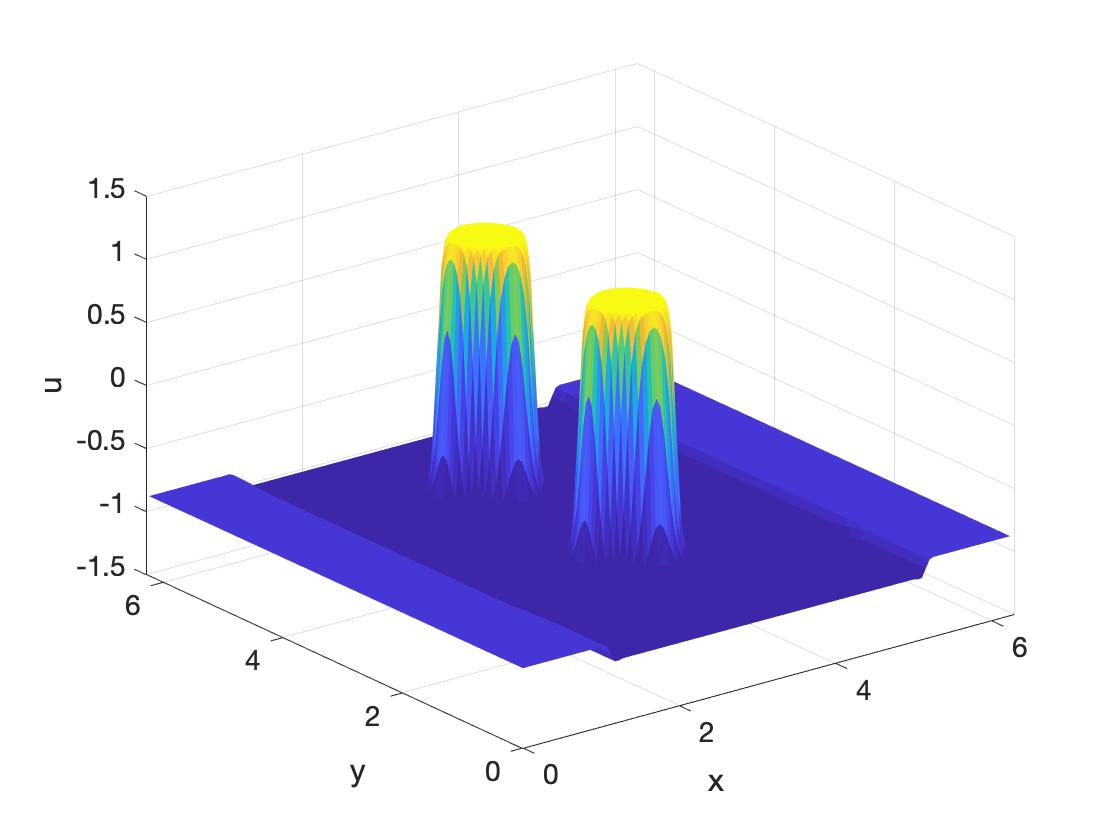}}
\subfigure[$u$ at $T=0.4$ \eqref{cutoff:allen}-\eqref{cutoff:allen2}.]{
\includegraphics[width=0.30\textwidth,clip==]{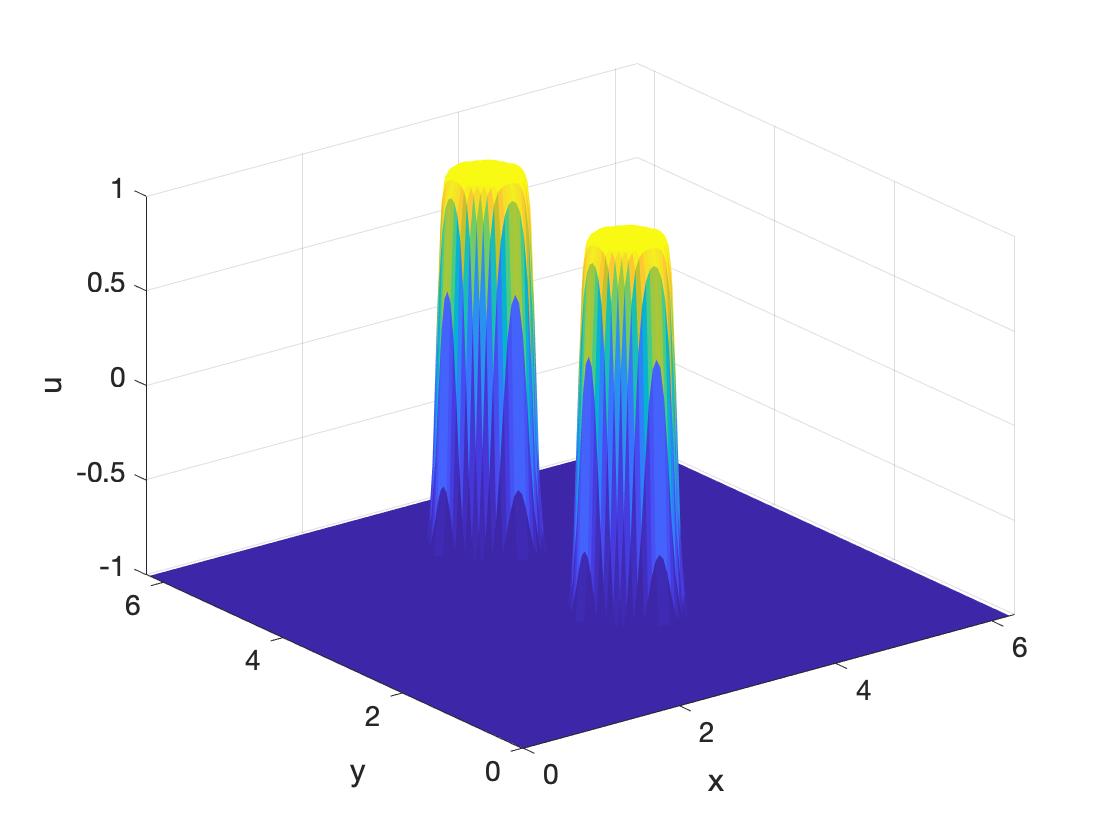}}
\subfigure[$\lambda$ at $T=0.4$.]{
\includegraphics[width=0.30\textwidth,clip==]{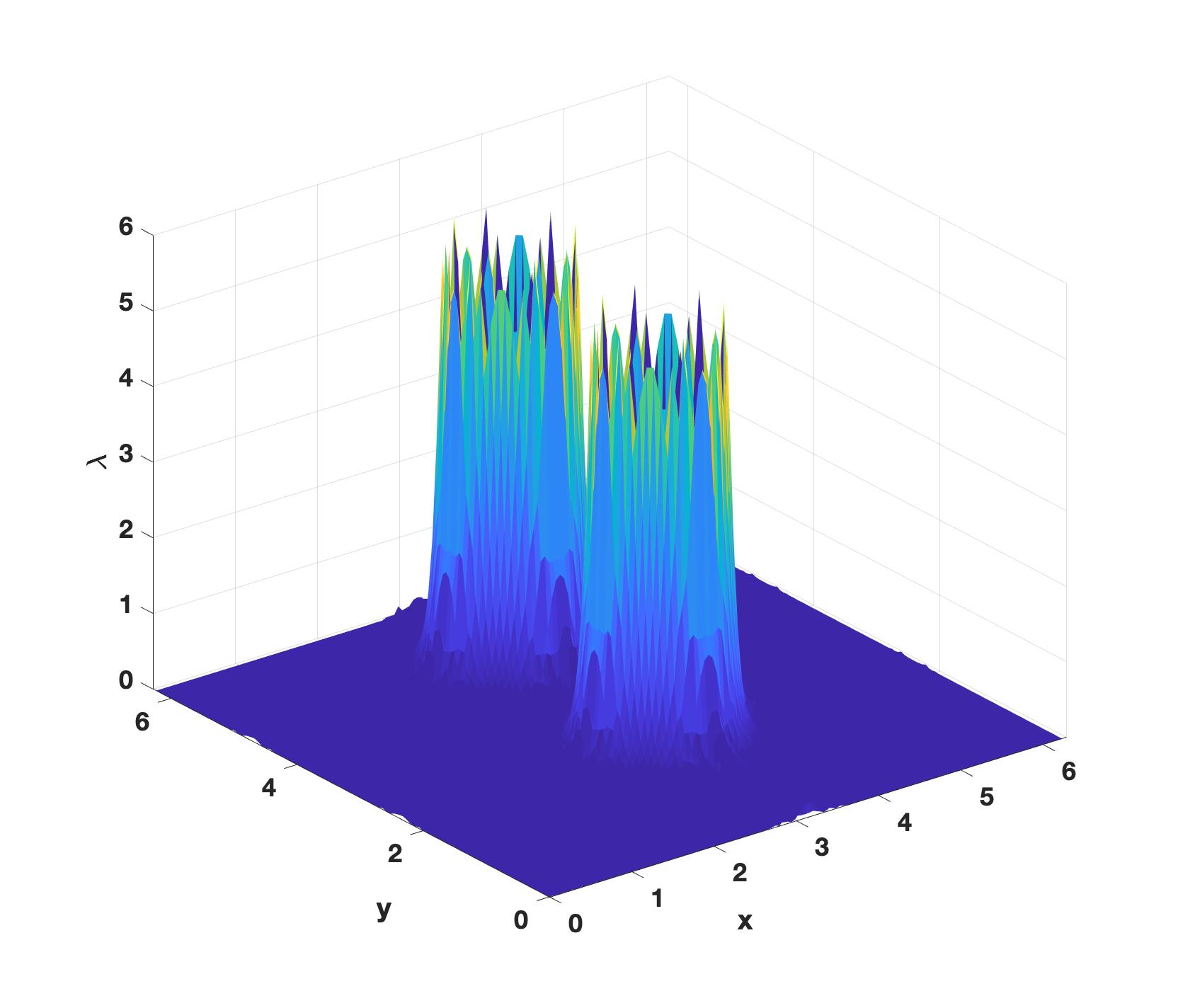}}
\caption{(a)-(d): Numerical solutions at $T=0.08,0.4$ computed by  \eqref{cutoff:allenB}.  (b)-(c) and (e)-(f): numerical solutions and Lagrange multiplier $\lambda$ at $T=0.08, 0.4$ computed by \eqref{cutoff:allen}-\eqref{cutoff:allen2}.}\label{allen:compare}
\end{figure}

\begin{figure}[htbp]
\centering
\includegraphics[width=0.45\textwidth,clip==]{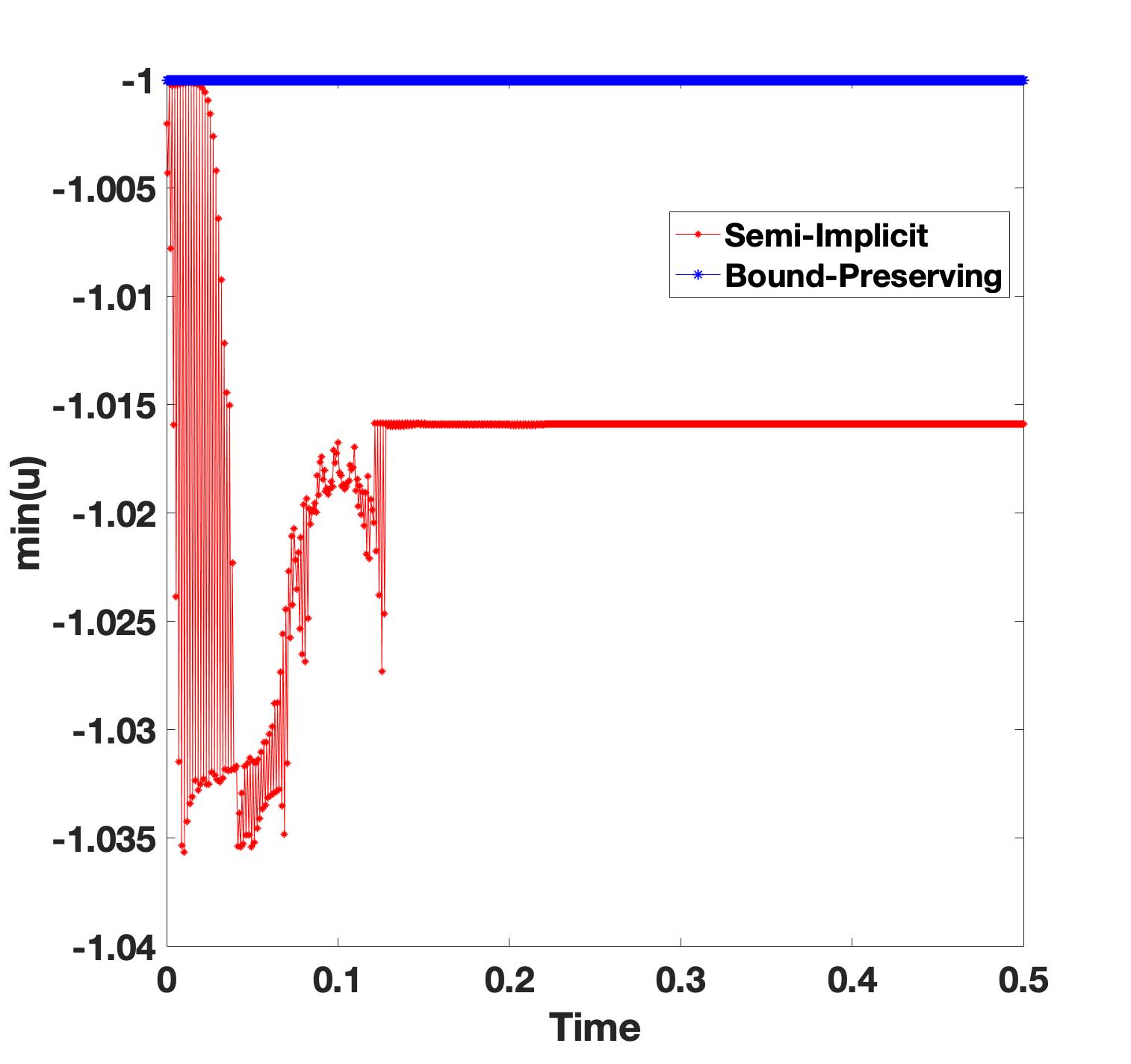}
\includegraphics[width=0.45\textwidth,clip==]{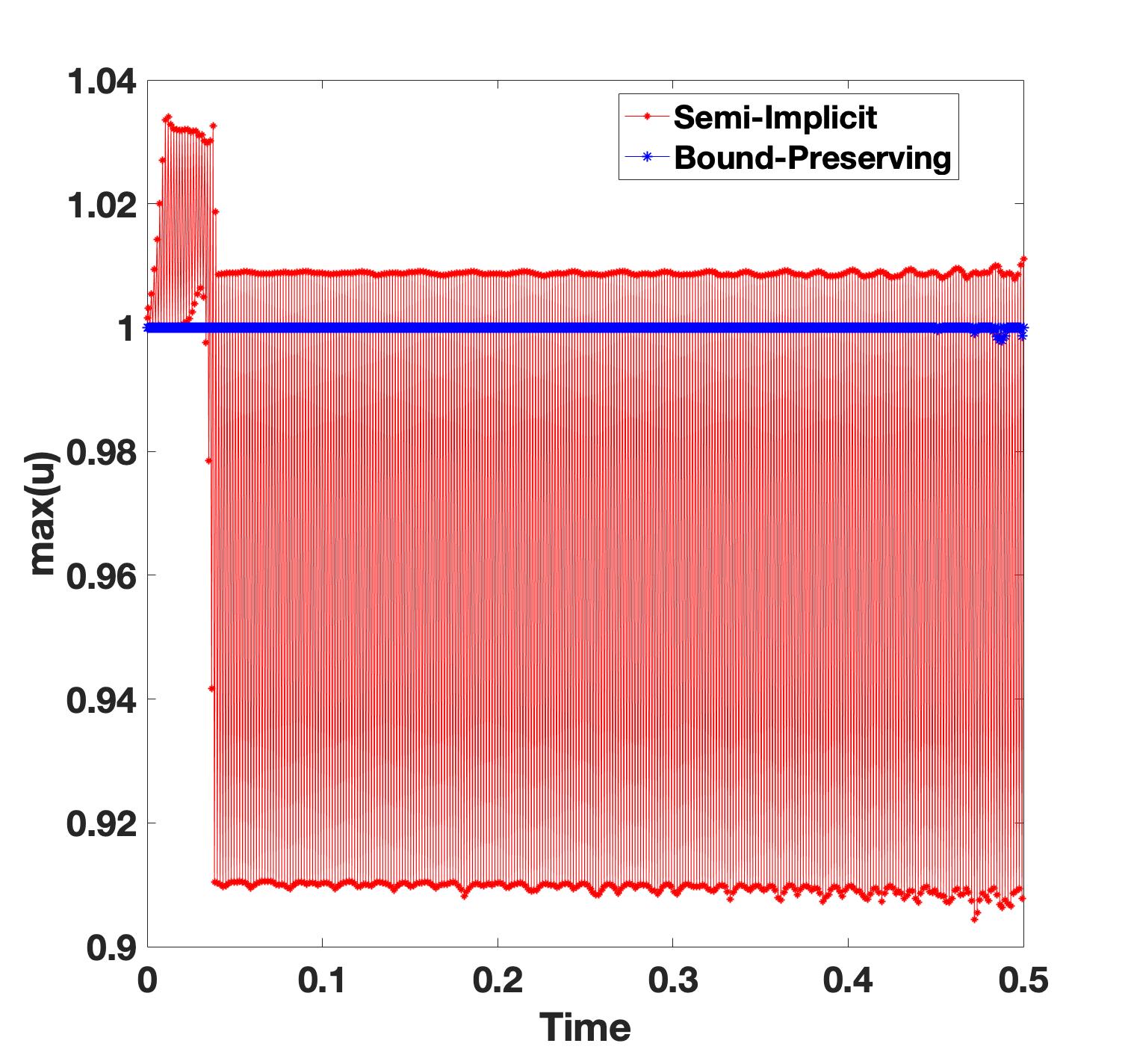}
\caption{Evolution of $\max\{u\}$ and $\min\{u\}$ with respect to time for the  semi-implicit scheme \eqref{cutoff:allenB} and the bound-preserving scheme  \eqref{cutoff:allen}-\eqref{cutoff:allen2}.}\label{max:bound}
\end{figure}

In Figure.\;\ref{allen:compare}, we plot the numerical solution  $u$ at $T=0.08$ and $T=0.4$ using the  semi-implicit scheme \eqref{cutoff:allenB} and the bound-preserving scheme  \eqref{cutoff:allen}-\eqref{cutoff:allen2}. It is observed that the numerical solution by the bound-preserving scheme  stays within $[-1,1]$, while that by  the  semi-implicit scheme \eqref{cutoff:allenB} violates this property. The  Lagrange multiplier $\lambda$ by  the bound-preserving scheme  \eqref{cutoff:allen}-\eqref{cutoff:allen2}  are also shown in  Figure.\;\ref{allen:compare}. 
 In Fig.\;\ref{max:bound}, we plot the evolution of $\max\{u\}$ and $\min\{u\}$ by both schemes.

\subsection{Cahn-Hilliard equation} 
\begin{figure}[htbp]
\centering
\includegraphics[width=0.50\textwidth,clip==]{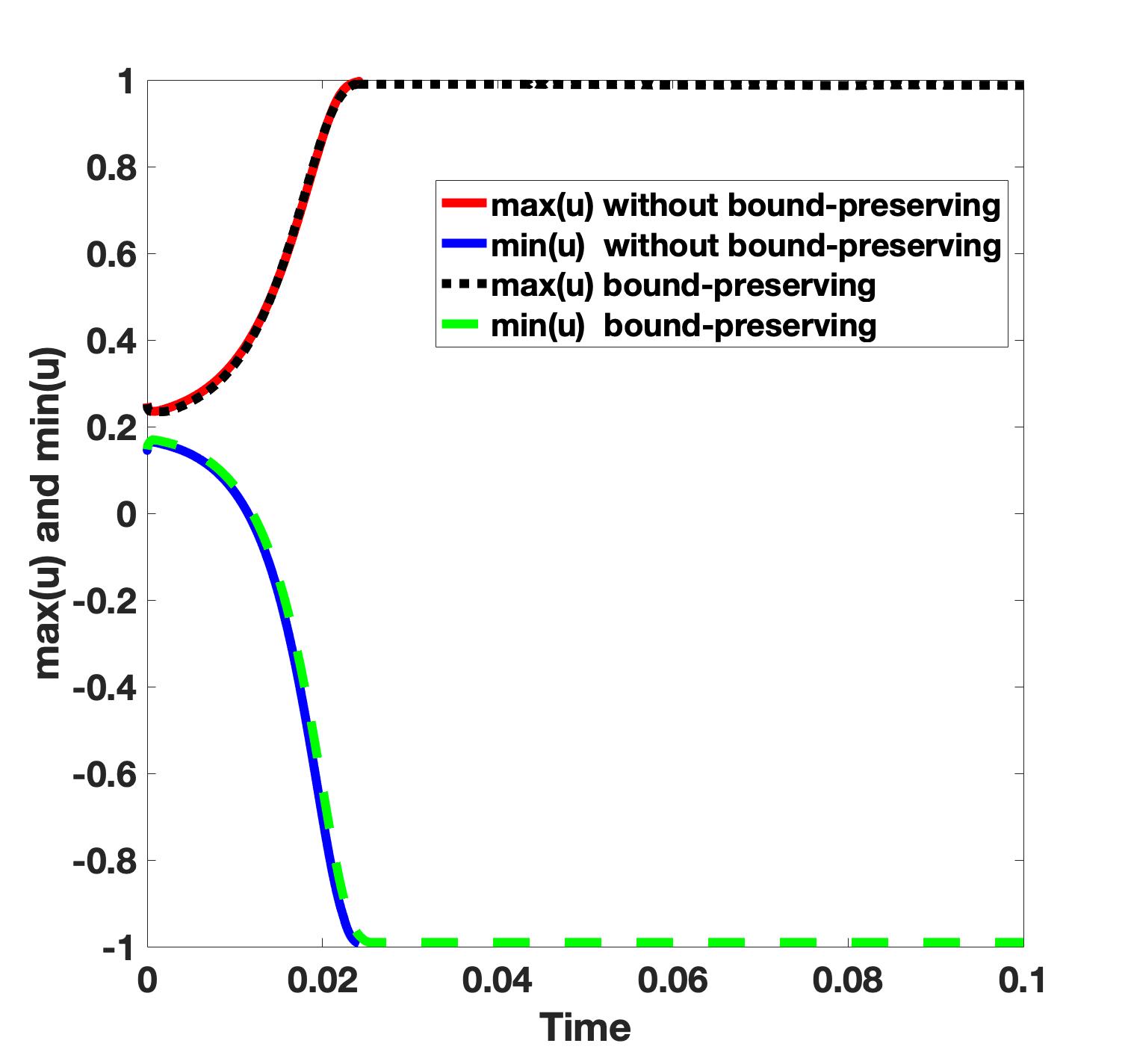}
\caption{The evolution of $\max_{i,j}{u^n_{i,j}}$ and $\min_{i,j}{u^n_{i,j}}$ with respect to time computed by the scheme  \eqref{bdf2:CH:LMB} up to $t\approx .025$, and  by the scheme \eqref{bdf2:CH:LM}-\eqref{bdf2:CH:LM2} up to $t=0.1$.}\label{bound:rho}
\end{figure}

\begin{figure}[htbp]
\centering
\subfigure[$t=0.001$]{
\includegraphics[width=0.23\textwidth,clip==]{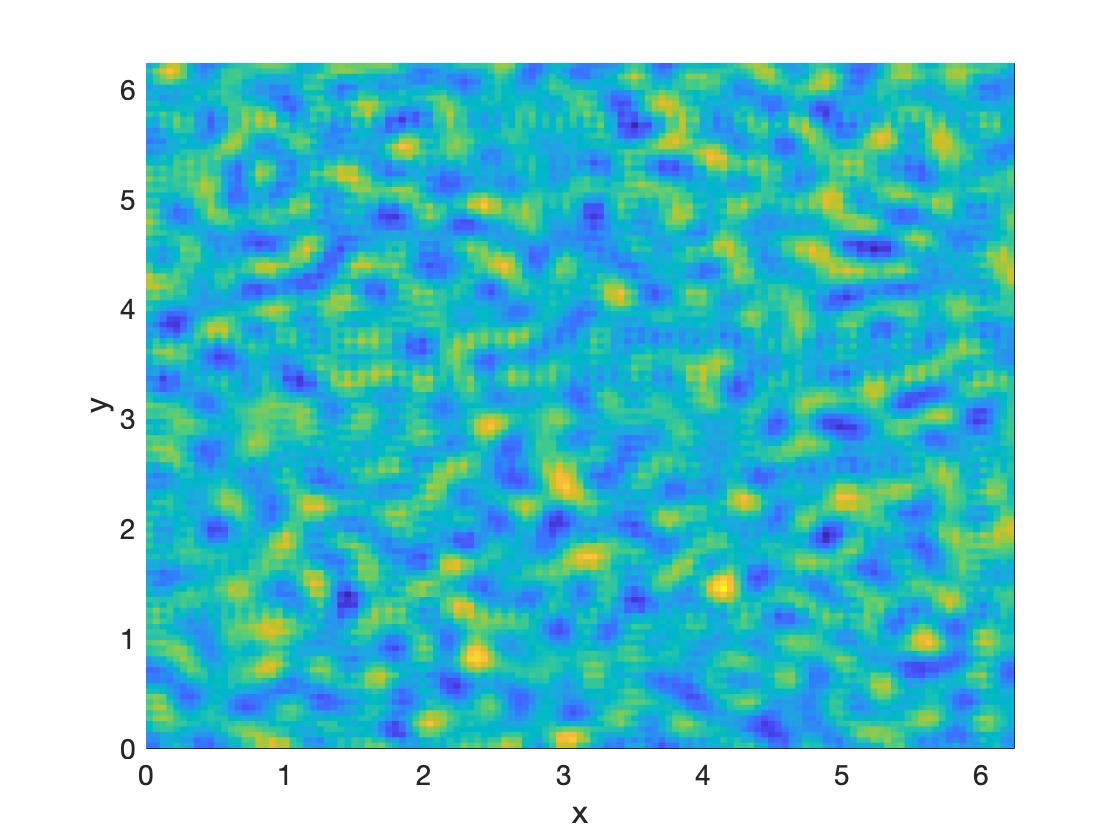}}
\subfigure[$t=0.02$]{
\includegraphics[width=0.23\textwidth,clip==]{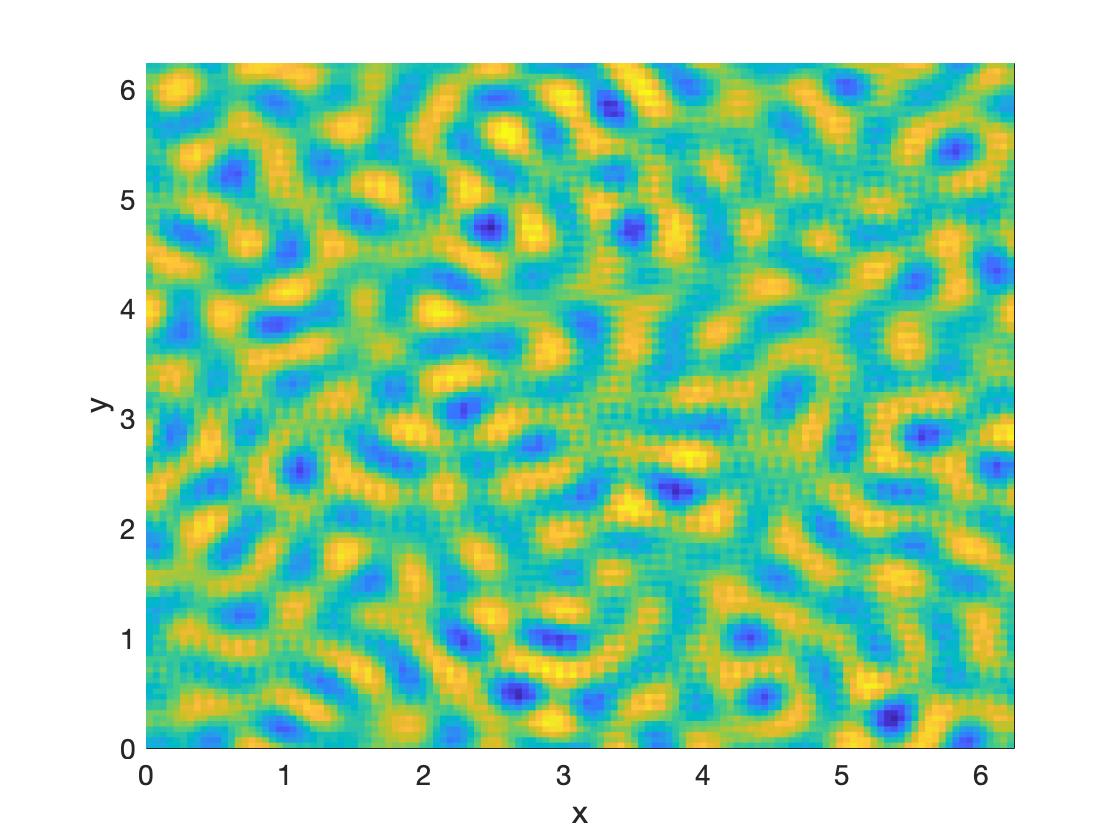}}
\subfigure[$t=0.05$]{
\includegraphics[width=0.23\textwidth,clip==]{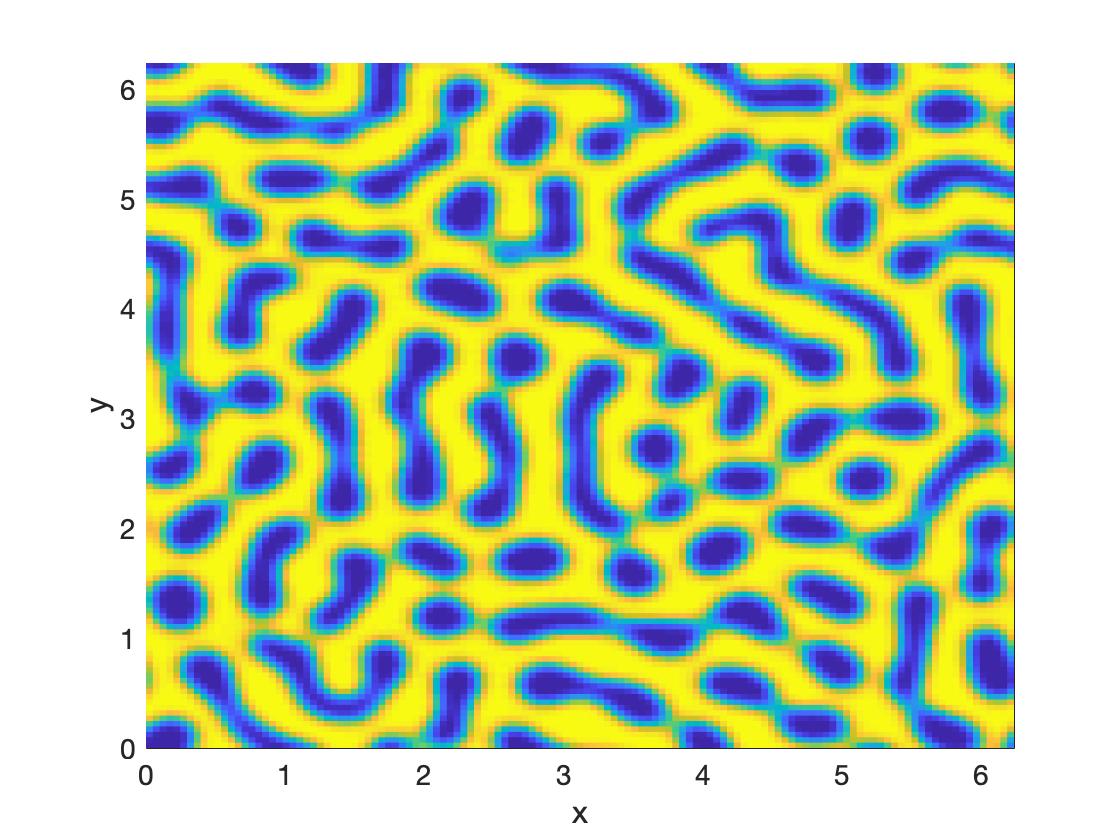}}
\subfigure[$t=0.1$]{
\includegraphics[width=0.23\textwidth,clip==]{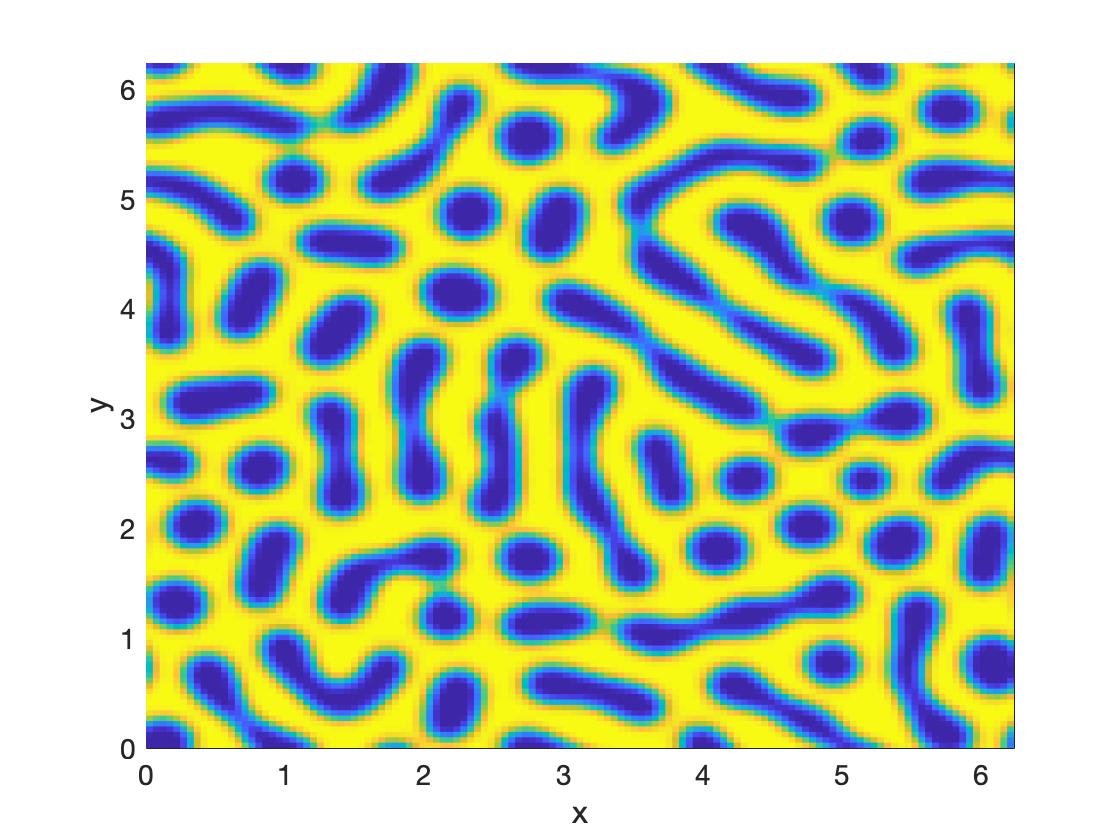}}
\caption{Numerical solutions of Cahn-Hilliard equation  at $t=0.001,0.02, 0.05, 0.1$ computed by  the scheme \eqref{bdf2:CH:LM}-\eqref{bdf2:CH:LM2}. }\label{coarse:CH}
\end{figure}

We now consider the Cahn-Hilliard equation \eqref{CH} with the initial condition 
\begin{equation}\label{ini_rand}
u_0(x,y)=0.2+0.05\,{\rm rand}(x,y),
\end{equation}
where  function ${\rm rand}(x,y)$ is a uniformed distributed random function with values in $(-1,1)$. We set $\theta_0=5$ and  $\eps=0.1$, and use    $\delta t=10^{-5}$  with $128^2$ Fourier modes in $(0,2\pi)^2$. 
We first use the following semi-implicit scheme 
\begin{equation}\label{bdf2:CH:LMB}
	\begin{split}
		& \frac{3{u}^{n+1}-4u^n+u^{n-1}}{2\delta t} =\Grad\cdot(M(2u^n-u^{n-1})\Grad\mu^{n+1}) ,\\
		&\mu^{n+1}= -\eps^2\Delta u^{n+1} +\ln(1+2u^n-u^{n-1})-\ln(1-2u^n+u^{n-1})-\theta_0(2u^n-u^{n-1}),
	\end{split}
\end{equation}
and found that it blows up  at $t\approx .025$ when $\|u^n\|_{l^\infty}>1$ due to the singular potential. We then use the bound preserving scheme \eqref{bdf2:CH:LM}-\eqref{bdf2:CH:LM2} with $\delta =0.01$ to compute up to $t=0.1$, and plot
in Fig.\;\ref{bound:rho}  the evolution of  $\max_{\bz\in\Sigma_N} u^{n}(\bz)$ and $\min_{\bz\in\Sigma_N} u^{n+1}(\bz)$ by the scheme \eqref{bdf2:CH:LMB} up to $t\approx .025$, and by the scheme \eqref{bdf2:CH:LM}-\eqref{bdf2:CH:LM2} up to $t=0.1$. 
 In Fig.\;\ref{coarse:CH}, we plot the numerical solutions  at various times which depict the coarsening process.

\subsection{Fokker-Planck equation}
As the final example, we consider
the Fokker-Planck equation \eqref{fockker-planck} with periodic boundary condition whose solution remains in $[0,1]$ and is mass preserving. We present below  simulations of  \eqref{fockker-planck}  on the domain $(-2\pi,2\pi)$  with
the initial condition  $u(x,0)=-e^{-\frac{(x-1)^2}{0.4}}$ 
using three second-order schemes: a usual semi-implicit scheme
\begin{equation}\label{scheme:fockker:planckB}
	\begin{split}
		&\frac{3u^{n+1}-4u^{n}+u^{n-1}}{2\delta t} = \partial_x(x(2u^n-u^{n-1})(1-2u^n+u^{n-1})+\partial_x u^{n+1}) ,
	\end{split}
\end{equation}
the bound-preserving scheme \eqref{scheme:fockker:planck}-\eqref{scheme:fockker:planck2} and  the mass conservative, bound-preserving scheme \eqref{mass:fockker:planck}-\eqref{mass:fockker:planck2}.

\begin{figure}[htbp]
 	\centering
 	\subfigure[By the bound-preserving scheme.]{
 		\includegraphics[width=0.45\textwidth,clip==]{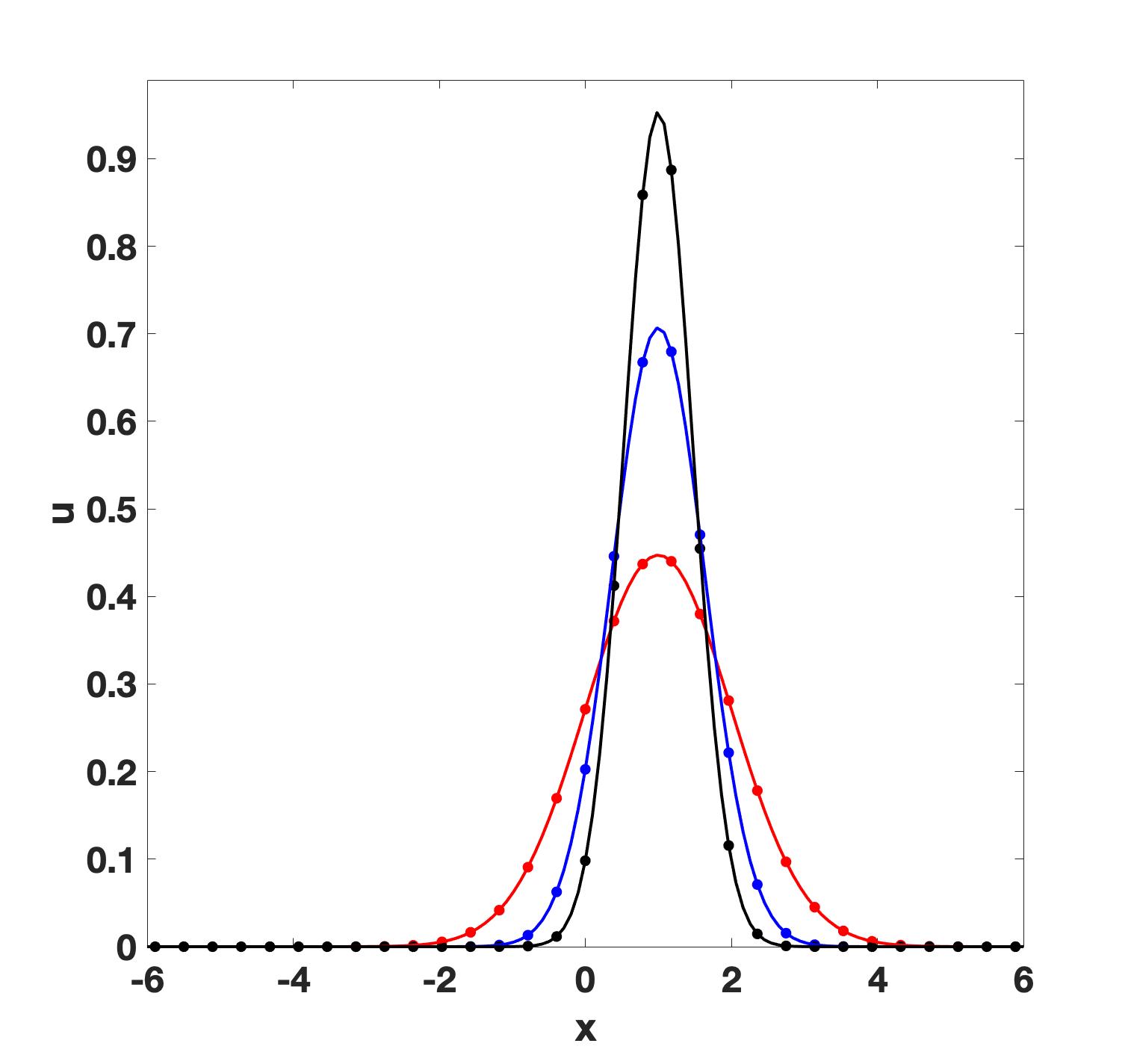}}
 	\subfigure[By the semi-implicit scheme.]{
 		\includegraphics[width=0.45\textwidth,clip==]{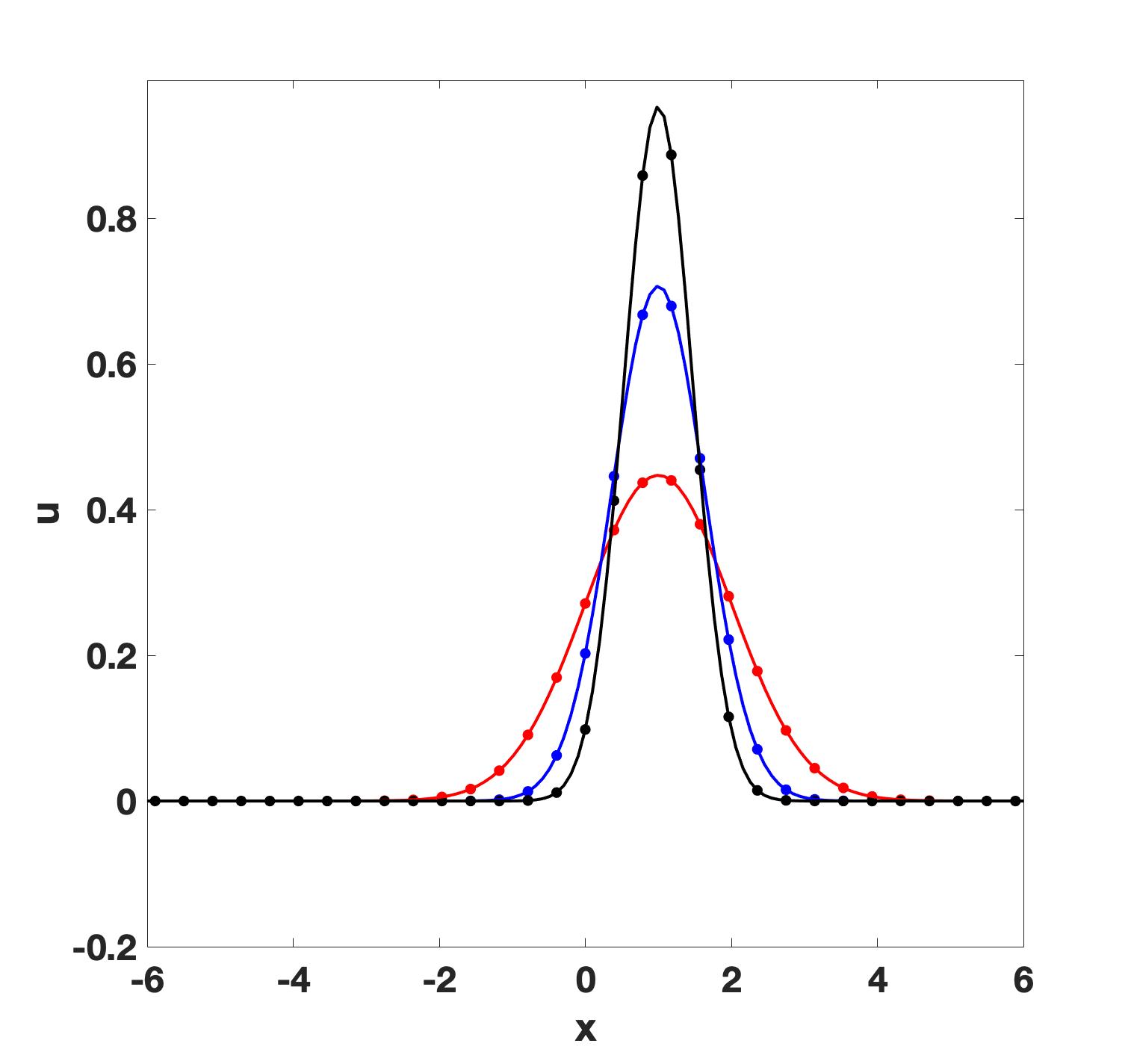}}
 	\subfigure[Evolutions of minimum values]{
 		\includegraphics[width=0.45\textwidth,clip==]{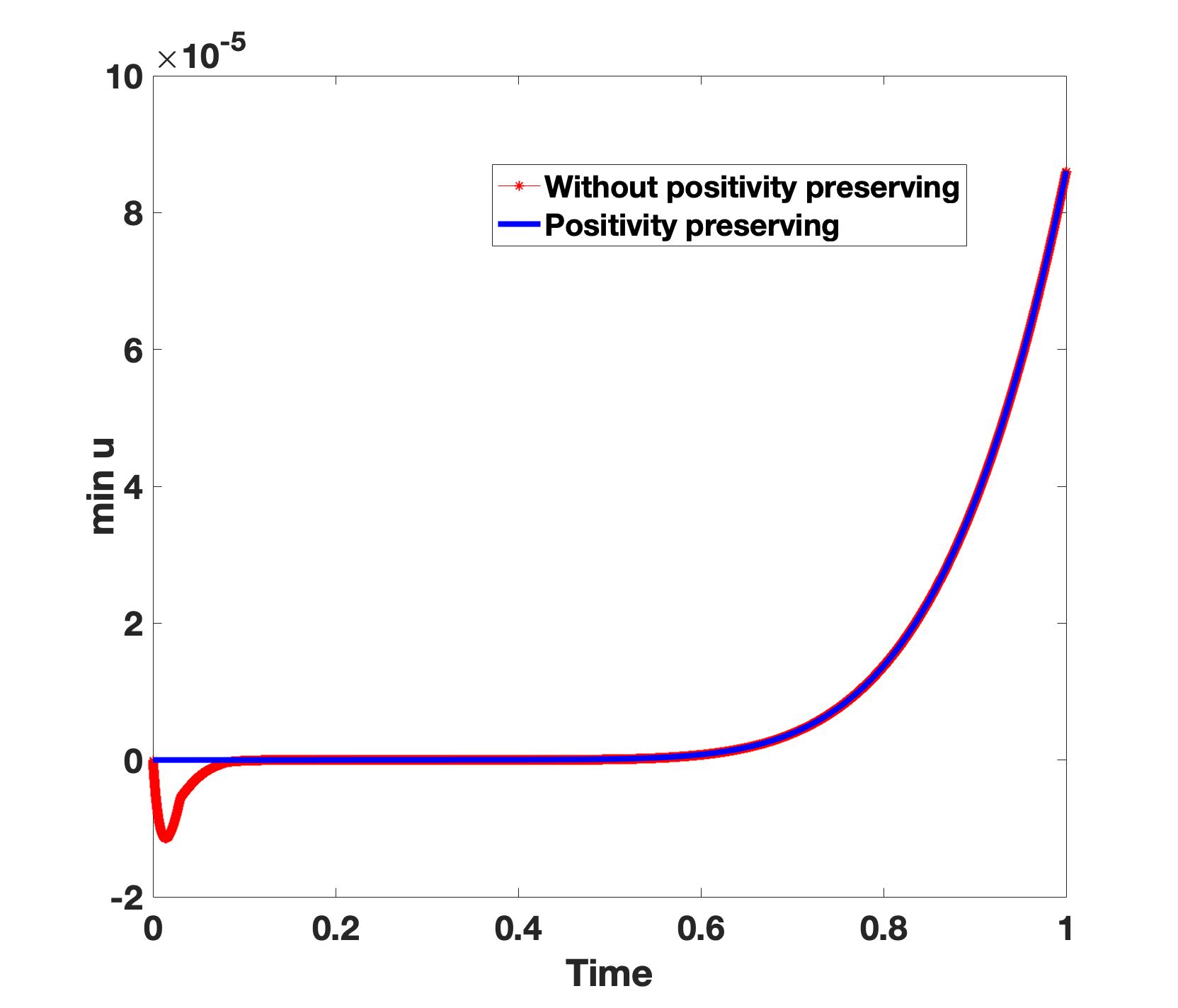}}
 	\subfigure[$\lambda$ at $t=0.01$]{
 		\includegraphics[width=0.42\textwidth,clip==]{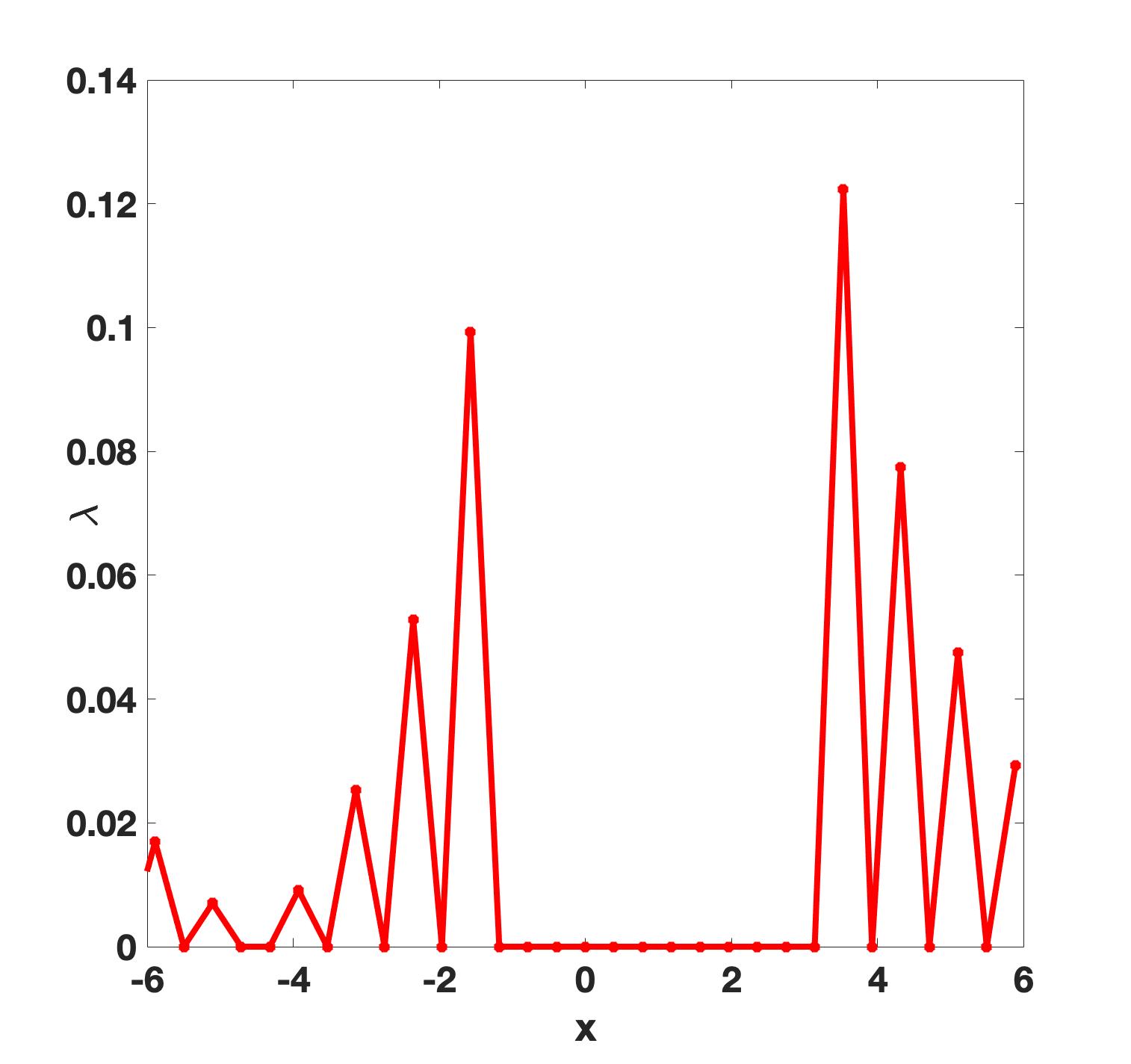}}
 	\caption{ (a)-(b) :Numerical solutions computed with $32$ Fourier modes plotted on the $256$ uniform grids  using   \eqref{scheme:fockker:planck}-\eqref{scheme:fockker:planck2}  and -\eqref{scheme:fockker:planckB}. (c): Evolutions of minimal values  using   \eqref{scheme:fockker:planck}-\eqref{scheme:fockker:planck2}  and -\eqref{scheme:fockker:planckB}. (d): Lagrange multiplier $lambda$ at $t=0.01$ using  \eqref{scheme:fockker:planck}-\eqref{scheme:fockker:planck2}.  }\label{bound:fock}
 \end{figure}

 In   Fig.\;\ref{bound:fock}, we plot the numerical results using the semi-implicit scheme \eqref{scheme:fockker:planckB}  and  the bound-preserving scheme \eqref{scheme:fockker:planck}-\eqref{scheme:fockker:planck2} with 32 Fourier modes and $\delta t=10^{-4}$. We observe that while the two numerical solutions look very similar, the minimum value by the semi-implicit scheme \eqref{scheme:fockker:planckB}  does become negative in a short period at the beginning, while   the numerical solutions by \eqref{scheme:fockker:planck}-\eqref{scheme:fockker:planck2} remain in $[0,1]$. 
 
\begin{figure}[htbp]
\centering
\subfigure[Evolution of mass]{
\includegraphics[width=0.40\textwidth,clip==]{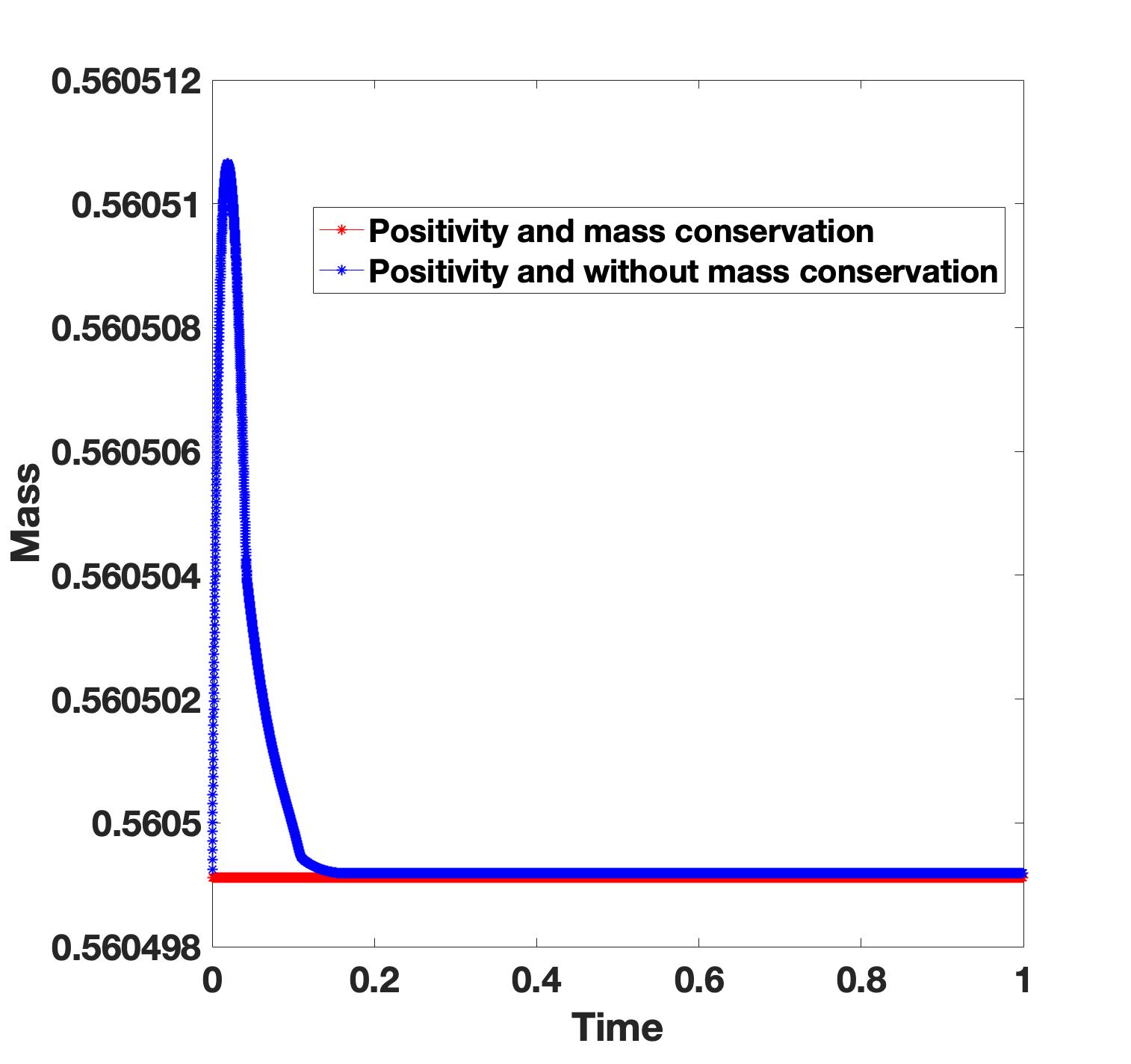}}
\subfigure[Solution profiles: $u$ at $t=0.01, 0.1, 0.4$.]{
\includegraphics[width=0.40\textwidth,clip==]{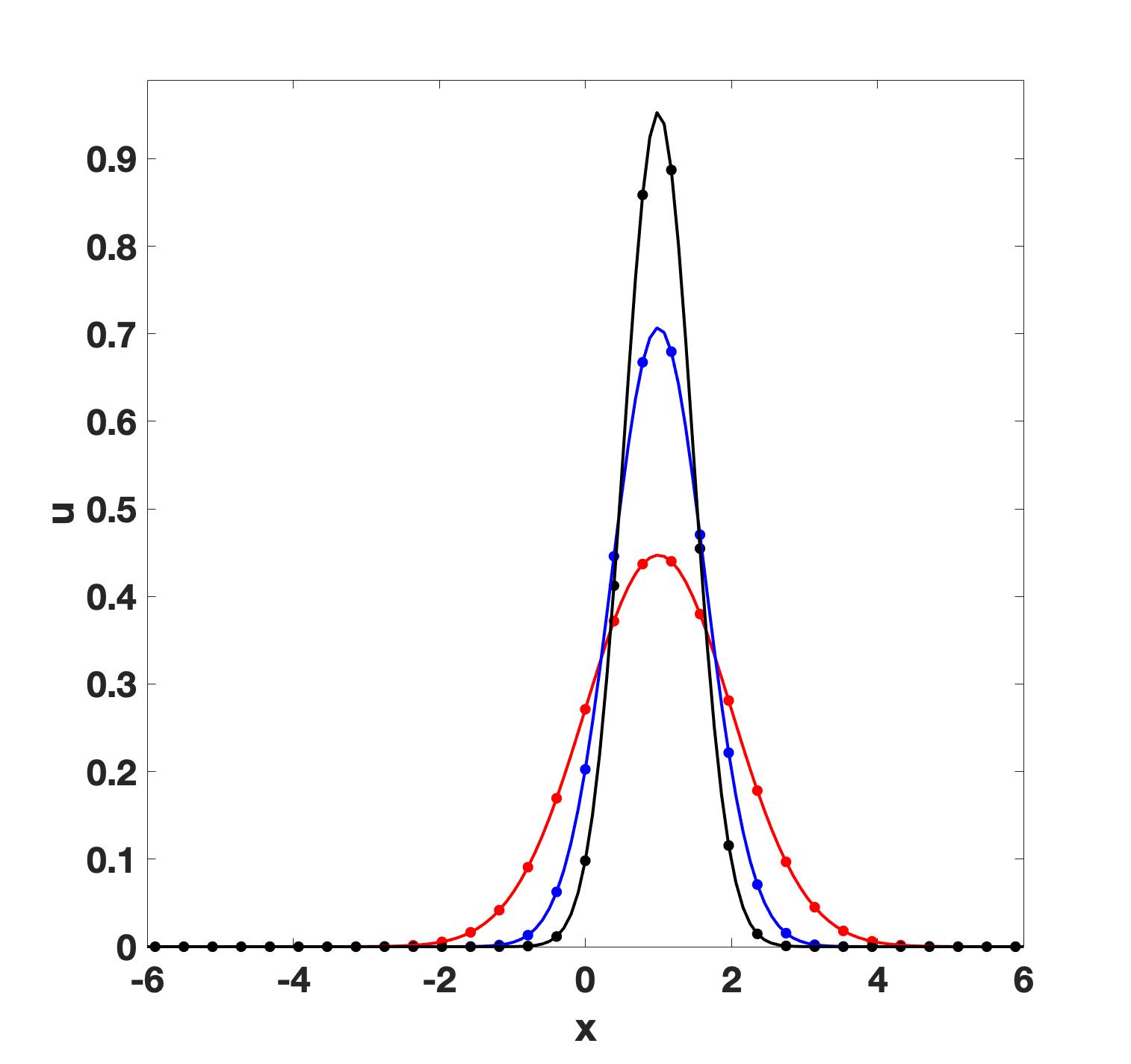}}
\subfigure[Iteration number ]{
\includegraphics[width=0.40\textwidth,clip==]{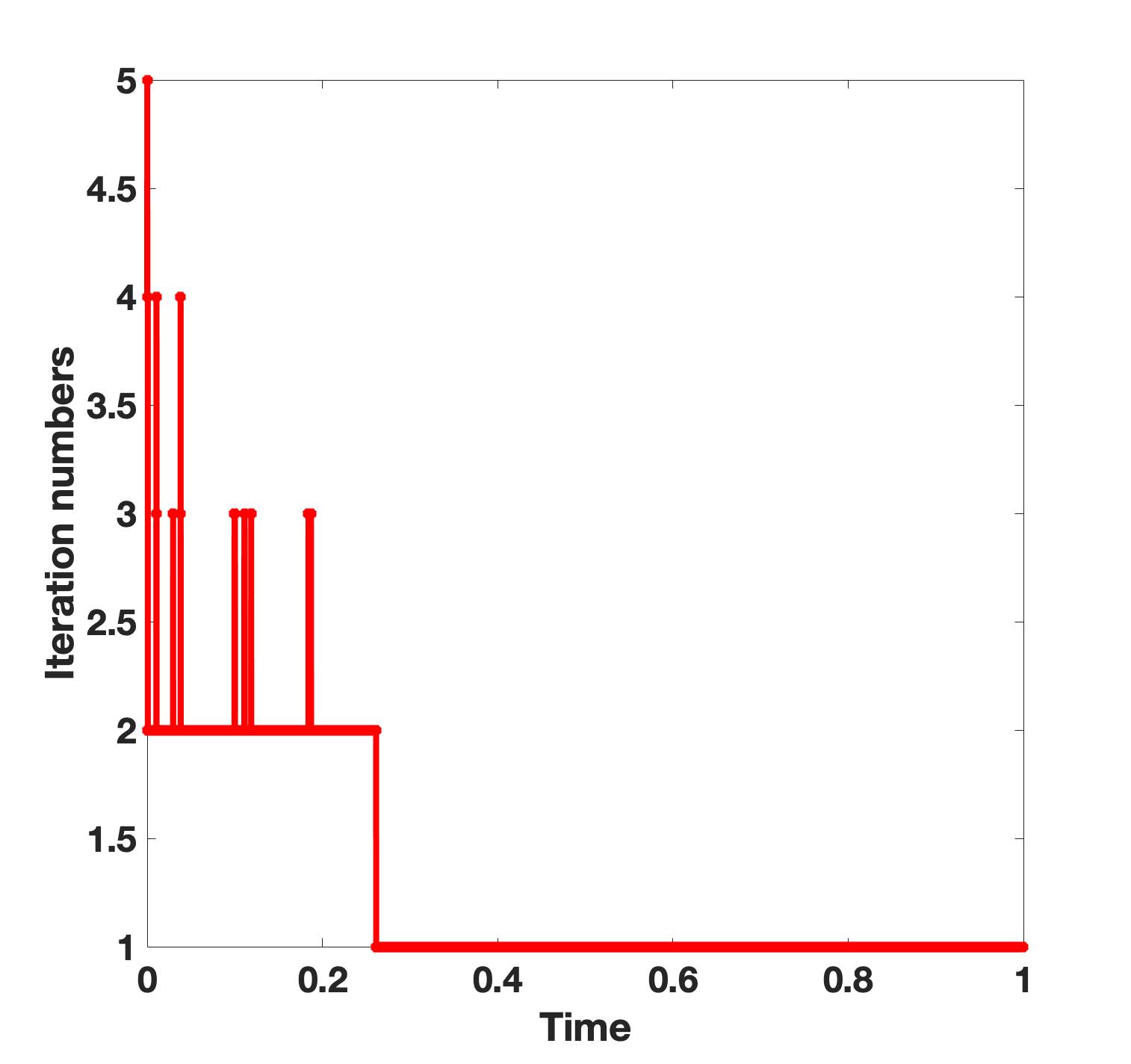}}
\caption{ (a): Evolution of mass  by  \eqref{scheme:fockker:planck}-\eqref{scheme:fockker:planck2} and  \eqref{mass:fockker:planck}-\eqref{mass:fockker:planck2}. (b): Solution profiles by \eqref{mass:fockker:planck}-\eqref{mass:fockker:planck2}. (c): Iteration numbers for solving $\xi^{n+1}$ at each time step of \eqref{mass:fockker:planck}-\eqref{mass:fockker:planck2}.} \label{mass:bound:fock}
\end{figure}

 In   Fig.\;\ref{mass:bound:fock}, we plot the numerical results using the bound-preserving scheme \eqref{scheme:fockker:planck}-\eqref{scheme:fockker:planck2}  and the  mass conservative, bound-preserving scheme \eqref{mass:fockker:planck}-\eqref{mass:fockker:planck2} with 32 Fourier modes and $\delta t=10^{-4}$. 
We observe that \eqref{scheme:fockker:planck}-\eqref{scheme:fockker:planck2} can not preserve mass, while   \eqref{mass:fockker:planck}-\eqref{mass:fockker:planck2}  preserves mass exactly.   Only a few iterations are needed to compute the  Lagrange multiplier $\xi$ at each time step by using the mass conservative, bound-preserving scheme \eqref{mass:fockker:planck}-\eqref{mass:fockker:planck2} .

\section{Concluding remarks}
We constructed efficient and accurate bound and/or mass preserving schemes for  a class of semi-linear and quasi-linear parabolic equations using the Lagrange multiplier approach.

First, we constructed a class of multistep IMEX schemes \eqref{high:bound:lag:1}-\eqref{high:bound:lag:2} for the semi-discrete problem \eqref{full:dis} with a Lagrange multiplier to enforce bound preserving, which is an approximation to the original PDE \eqref{strong}.  Hence, the scheme  \eqref{high:bound:lag:1}- \eqref{high:bound:lag:2}  is a $k$-th order approximation in time for both  \eqref{full:dis} and  \eqref{strong}.
In particular, the  \eqref{high:bound:lag:1}-\eqref{high:bound:lag:2} can be very useful if one is interested in the discrete problem  \eqref{full:dis} without a background PDE.

Then, we pointed out in {\bf Remark} 2.1 that by dropping out the term $B_{k-1}(\lambda_h^ng'(u_h^n))$ in \eqref{high:bound:lag:1} and \eqref{high:bound:lag:2}, 
we recover  the usual cut-off scheme which is a $k$-th order approximation in time for   \eqref{strong}, but only a first-order approximation in time for  \eqref{full:dis}. 
Thus, our presentation provided an alternative interpretation of the cur-off approach, and moreover, allowed us to construct  new cut-off  implicit-explicit (IMEX) schemes with mass conservation.

 We also established  some  stability results involving norms with derivatives under a general setting, and derived  optimal error estimates for a second-order bound preserving scheme with a hybrid spectral discretization in space. 
 
 Finally, we applied our approach to several typical PDEs which preserve bound and/or mass, and presented ample numerical results to validate our approach. The approach presented in this paper is quite general and can be used to develop bound preserving schemes for other bound preserving PDEs such as the Keller-Segel equations \cite{Kel.S70}.

\bibliographystyle{siamplain}
 \bibliographystyle{plain}
\bibliography{references}
\end{document}